\documentclass[5p]{elsarticle}
\usepackage{hyperref}
\usepackage{natbib}
\usepackage[english]{babel}
\usepackage{amssymb}
\usepackage{amsmath} 
\usepackage{framed} 
\usepackage{cancel}

\usepackage{multicol} 
\usepackage{nomencl} 
\renewcommand*\nompreamble{\begin{multicols}{2}}
\renewcommand*\nompostamble{\end{multicols}}
\usepackage{etoolbox}

\renewcommand\nomgroup[1]{%
  \item[\itshape 
  \ifstrequal{#1}{A}{Variables}{%
  \ifstrequal{#1}{B}{Parameters}{ 
  \ifstrequal{#1}{C}{Subscripts and sets}{
  \ifstrequal{#1}{D}{Abbreviations}{ }}}}%
]\vspace{10pt}} 

\makeatletter
\providecommand{\doi}[1]{%
  \begingroup
    \let\bibinfo\@secondoftwo
    \urlstyle{rm}%
    \href{http://dx.doi.org/#1}{%
      doi:\discretionary{}{}{}%
      \nolinkurl{#1}%
    }%
  \endgroup
}
\makeatother

\addto\captionsenglish{}

\journal{arXiv.org}

\begin{document}

\begin{frontmatter}

\title{Impact of different time series aggregation methods on optimal energy system design}

\author[addFZJ]{Leander Kotzur\corref{corLK}}
\author[addFZJ]{Peter Markewitz}
\author[addFZJ]{Martin Robinius}
\author[addFZJ,addCFC]{Detlef Stolten}

\cortext[corLK]{Corresponding author. Email: l.kotzur@fz-juelich.de}
\address[addFZJ]{ Institute of Electrochemical Process Engineering (IEK-3), Forschungszentrum J\"ulich GmbH, Wilhelm-Johnen-Str., 52428 J\"ulich, Germany}
\address[addCFC]{Chair for Fuel Cells, RWTH Aachen University, c/o Institute of Electrochemical Process Engineering (IEK-3), Forschungszentrum J\"ulich GmbH, Wilhelm-Johnen-Str., 52428 J\"ulich, Germany}

\begin{abstract}
Modelling renewable energy systems is a computationally-demanding task due to the high fluctuation of supply and demand time series. To reduce the scale of these, this paper discusses different methods for their aggregation into typical periods. Each aggregation method is applied to a different type of energy system model, making the methods fairly incomparable.

To overcome this, the different aggregation methods are first extended so that they can be applied to all types of multidimensional time series and then compared by applying them to different energy system configurations and analyzing their impact on the cost optimal design. 

It was found that regardless of the method, time series aggregation allows for significantly reduced computational resources. Nevertheless, averaged values lead to underestimation of the real system cost in comparison to the use of representative periods from the original time series. The aggregation method itself – e.g., k-means clustering – plays a minor role. More significant is the system considered: Energy systems utilizing centralized resources require fewer typical periods for a feasible system design in comparison to systems with a higher share of renewable feed-in. Furthermore, for energy systems based on seasonal storage, currently existing models’ integration of typical periods is not suitable.

\end{abstract}


\end{frontmatter}

\begin{table*}[!t]   
\begin{framed}
\printnomenclature
\end{framed}
\end{table*}

%
\section{Introduction}
\label{sec:Introduction}
Developing an energy system design that minimizes costs and environmental impact is a complex task due to the spatial and temporal gap between energy production and demand. In consequence, optimization algorithms are required for  solving these design problems \cite{Banos2011, Petruschke2014, Stadler2014, Kwon2016, Samsatli2015, Merkel2015, Milan2012, Lauinger2016, Merei2016, Wang2015,Palzer2014}. 

However, the algorithms used hitherto are computationally demanding: The size of the input data directly influences that of the related optimization problem, and with it the requirement for processing resources. For this reason, it is often necessary to simplify the design problem in advance. 



Therefore, different options for complexity reduction exist and include: Spatial aggregation which reduces the number of nodes in an energy system network \cite{Mancarella2014}; simplifying the technology models by reducing nonlinearities or discontinuities so as to avoid non-convexity of the program \cite{Geidl2007, Milan2015}; and temporal aggregation, which creates typical periods representing the original input time series. 

\begin{figure}[h]
	\centering0
  \includegraphics{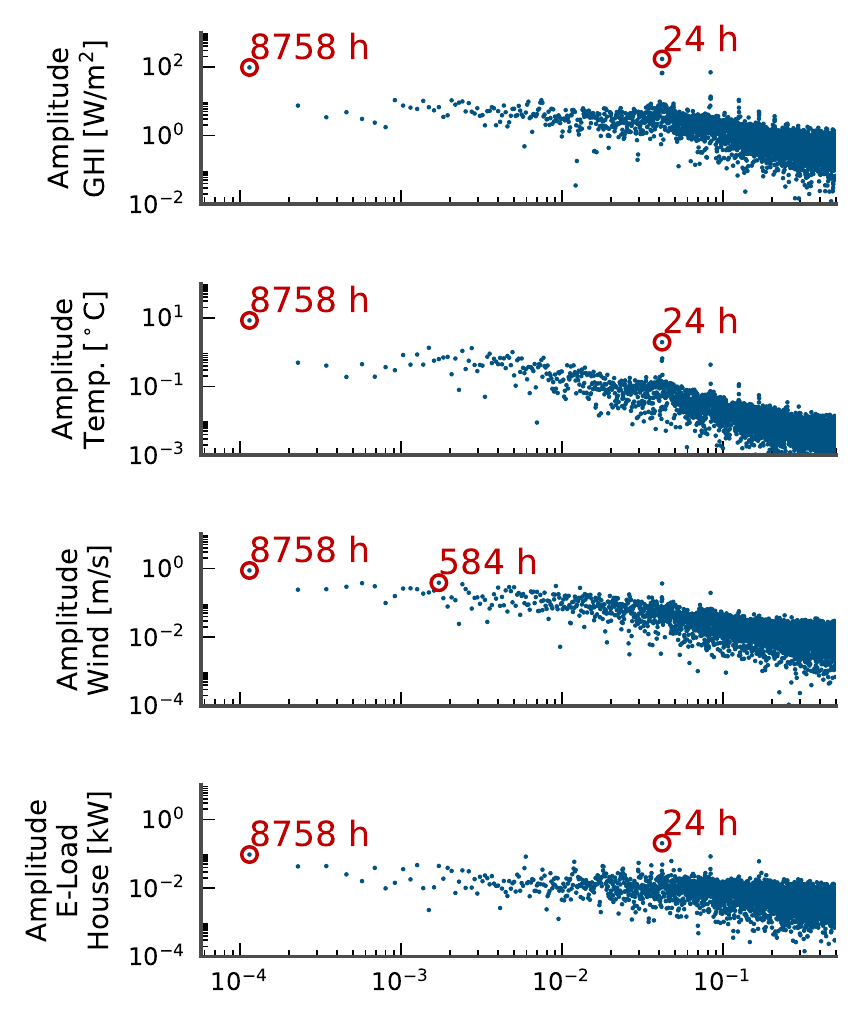}
	\caption{Fast Fourier Transformation of the Global Horizontal Irradiance (GHI) \nomenclature[D]{GHI}{Global Horizontal Irradiance}{}{}, the temperature and the wind speed of a test reference year (Location: Bad Marienberg, Germany) \cite{DWD2012} and a representative electrical load profile of a residential building (Profile 1) \cite{Tjaden2015} }
	\label{fig:FFT}
\end{figure}

The creation of recurring periods is popular because of the existing patterns in the hourly, daily and seasonal variation for the majority of design relevant time series. Therefore, it is reasonable to reduce redundant data until the minimal representative data set required for the problem is reached. Figure \ref{fig:FFT} visualizes this redundancy by showing the result of a Fast Fourier Transformation (FFT) of different time series that are typically required for an energy system design. The frequencies with the highest amplitudes are highlighted and are, as anticipated, the daily and annual variations.

For this reason, many different methods for the selection of typical periods have been presented. Aside from custom exact optimization methods  \cite{Golling2012,Poncelet2016}, and graphical methods \citep{Ortiga2011}, the majority use heuristic methods or greedy clustering algorithms for the aggregation of typical periods. Creating representative days by averaging time series, for example over a type of day defined by month or weekday, has been popular \cite{Mavrotas2008,Mehleri2013,Casisi2009,Lozano2009}. \cite{Lythcke2016} refers to it as time-chronological \textit{averaging}. Recent attempts use the \textit{k-means} clustering \cite{Adhau2014,Green2014,Bahl2016,Fazlollahi2014a}, \textit{hierarchical} clustering \cite{Merrick2016,Nahmmacher2016}, or \textit{k-medoids} clustering - either based on a greedy algorithm \cite{Rager2015,Stadler2016} or an exact solution of a MILP  \cite{Dominguez2011,Schuetz2017} - for the selection of typical periods. Nevertheless, each method is applied to a different system and it is difficult to identify which is the most suitable.
While \citet{Schuetz2016} compare methods for  building energy systems, they note that future research should focus on the appropriateness of clustering algorithms for different design applications. Moreover, the period length should be varied so as to assess the impact of storage effects.

A further difficulty is that the system considered determines the minimal required dataset. For renewable energy systems, a higher resolution of the input times-series is required than with their fossil counterparts \cite{Poncelet2014}. For conventional system design, it could be sufficient to reduce the dataset to a few time steps \cite{Bahl2016}, while for a storage based system design different typical weeks are required \cite{Harb2015,Nahmmacher2016}.



In summary, the following open research questions present themselves: 
\begin{itemize}
\item Which time-series aggregation method is best suited for which energy system design application?
\item What is the minimum number of aggregated time steps to model such a system?
\item What is an appropriate  period length - typical days or typical weeks?
\end{itemize}

To answer these, this paper is structured as follows: First different deterministic methods including \textit{k-mean} clustering, \textit{k-medoids} clustering and \textit{hierarchical} clustering as aggregation method are presented in section \ref{sec:02_Methods}, where the possibilities of adding extreme periods are also discussed. 
In section \ref{sec:03_Aggregation}, the aggregation methods are used to select four typical days of different time series that could be relevant for an energy system design. The aggregated profiles are then graphically analyzed and through accuracy indicators.
In the following, the different methods are applied in section \ref{sec:examples} to three design optimization problems of a heat and electricity supply system:
\begin{enumerate}
\item A cogeneration unit with a heat storage as benchmark system
\item A residential system based, amongst other elements, on photovoltaics and a heatpump
\item An island system with a high share of renewables with the support of different storage technologies
\end{enumerate}
To validate the methods, the results for different numbers of typical periods are compared to the optimal solution of the original optimization problems with the full time series and analysed in terms of their accuracy and computational load. 
Sections \ref{sec:Summary}, \ref{sec:Conclusion} and \ref{sec:Outlook} summarize, draw the principal conclusions and give an outlook for further research questions.

All methods introduced are published in the Python package \href{https://github.com/FZJ-IEK3-VSA/tsam/}{tsam - Time Series Aggregation Module} and can be easily applied and extended.

\section{Time series aggregation methods}
\label{sec:02_Methods}
The aim of time series aggregation is to merge a set of periods into groups such that the group members - the original periods - are as similar as possible. The group is then represented by a single period. 
The grouping of time-series is based in the most of the methods on a distance measure of the attributes between each group member.  For an accurate grouping, the raw input data must first be pre-processed into the right format (Section \ref{sec:Pre-processing}). On the basis of this, an aggregation method is then applied to create the groups (Section \ref{sec:aggregation-methods}). In terms of achieving a feasible system design, different variants of integrating extreme periods can be included afterwards (Section \ref{sec:extreme-periods}). Finally, the aggregated time series must be scaled back such that their average values fit the average values of the original time series (Section \ref{sec:scaling}). 

\begin{figure}[h]
	\centering
  \includegraphics[width=0.49\textwidth]{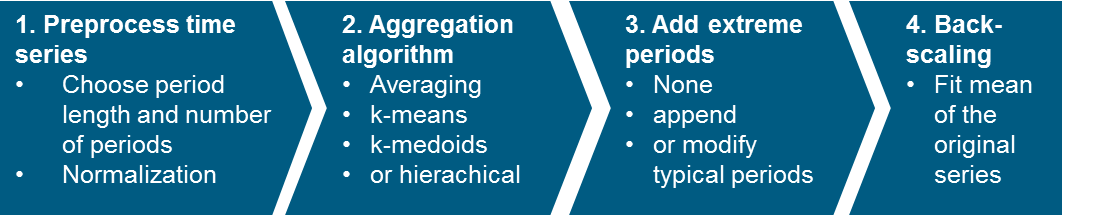}
	\caption{Different steps in the procedure of time series aggregation.}
	\label{fig:GeneralProcedure}
\end{figure}

The general procedure for time series aggregation applied in this paper is shown in Figure \ref{fig:GeneralProcedure}.

\subsection{Pre-processing the time-series}
\label{sec:Pre-processing}
\nomenclature[B]{$N$}{Size of an index set}{}{}
First, the input data is normalized in terms to evaluate all time series on the same scale. Each time series $\mathbf{x}_{a}^{'}$\nomenclature[B]{$x$}{Normalized candidate value}{}{} represents an attribute $a \in \{ 1,...,N_a\}$ \nomenclature[C]{$a$}{Attribute represented by a time series}{}{} e.g., an electrical demand profile or the measured solar horizontal irradiance at a certain location. The time series itself consists of raw data points $x_{a,t}^{'}$ where  $t \in \{ 1,...,N_t\}$\nomenclature[C]{$t$}{Time step index of the full series}{}{} constitutes a single time step. Different possibilities for normalization and standartization are presented by \citet{Rager2015}. In this work, the normalized time series $\mathbf{x}_{a}$ are calculated as follow
\begin{equation}
\mathbf{x}_{a} = \frac{\mathbf{x}_{a}^{'} - \min \mathbf{x}_{a}^{'}  }{ \max  \mathbf{x}_{a}^{'} - \min{  \mathbf{x}_{a}^{'} }} \quad   \forall \quad a \in \{ 1,...,N_a\}, 
\end{equation}
 which results in time series on the same scale $x_{a,t} \in \left[0,1  \right] $.

\begin{figure}[h]
	\centering
  \includegraphics[width=0.49\textwidth]{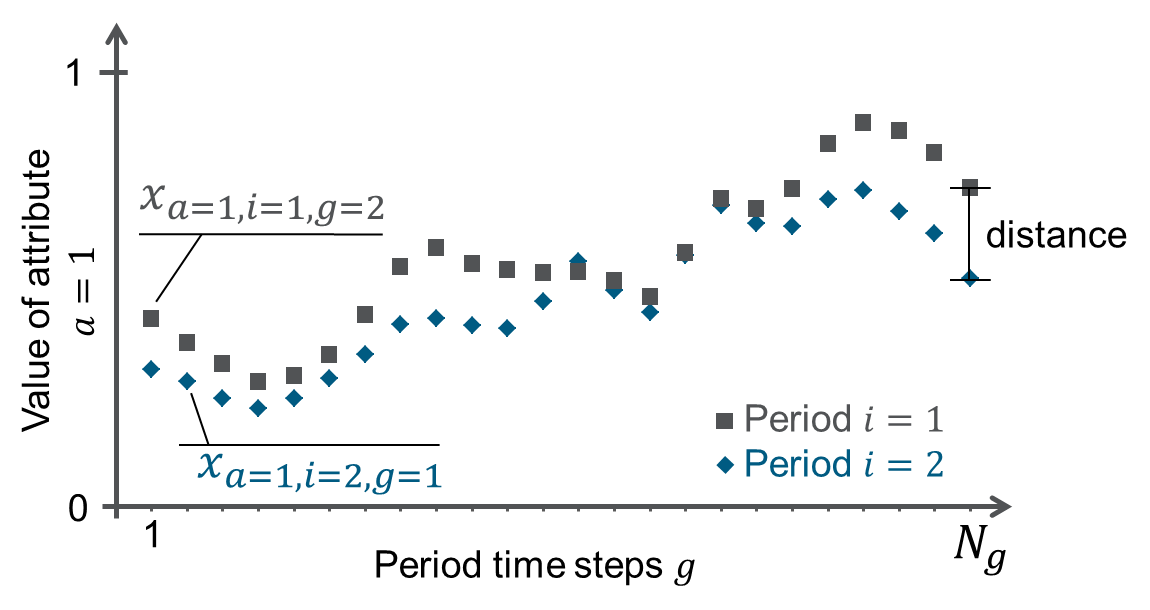}
	\caption{Illustration of the nomenclature for a time series $a=1$ with two candidate periods $i=1,2$ for a whole period.}
	\label{fig:PeriodsNomenclature}
\end{figure}

For the aggregation of typical periods, the scaled time-series are arranged into the candidate periods considered $i \in \{ 1,...,N_i \}$ \nomenclature[C]{$i$}{Candidate period index}{}{}, each consisting of the same number of time steps $g \in \{1,...,N_g \}$ \nomenclature[C]{$g$}{Time step index inside a period}{}{} with $N_a$ attributes. The nomenclature is illustrated in Figure \ref{fig:PeriodsNomenclature} where the normalized time-series of attribute $a=1$ is shown for two periods $i={1,2}$, which were originally successive. 

This reordering results in a matrix $L$ 
in which the number of columns is defined by the multiple of the number of period time steps $N_g$ and number of attributes $N_a$. The number of rows corresponds to the number of periods $N_i$: 
\begin{equation}
\label{eq:attribute-matrix}
\hspace*{-1cm}
L = \left[ \begin{array}{cccccc}
x_{1,1,1} & \cdots  & x_{1,N_g,1} & x_{2,1,1} & \cdots & x_{N_a,N_g,1} \\
\vdots & \ddots & \vdots & \vdots & \ddots & \vdots \\
x_{1,1,N_i} & \cdots  & x_{1,N_g,N_i} & x_{2,1,N_i} & \cdots & x_{N_a,N_g,N_i} \\
\end{array} \right] 
\end{equation}

A single row, described by vector $\pmb{x}_i$ represents a candidate period, also referred to as an observation point. 

For example: Pre-processing two hourly time series over a year  (${N_t = 8760}$), like the mentioned electrical demand profile and the measured solar horizontal irradiance (${N_a = 2}$),  for the application of a typical day approach (${N_g=24}$) would yield 365 days (${N_i = 365}$) and a $365$-by-$48$ matrix.

In the case that the raw time series length $N_t$ is no integer multiple of the period length $N_g$ and the number of periods $N_i$, the full time-series $x_{t,a}$ must be cut off or extended, such that its length becomes an integer multiple of the number of steps inside a candidate period $N_g$ and the number of periods $N_i$.

\subsection{Time-series aggregation methods}
\label{sec:aggregation-methods}

Based on the matrix introduced in equation \ref{eq:attribute-matrix}, different aggregation methods can be applied to group the $N_i$ independent candidate periods into clusters defined as $C_k$\nomenclature[B]{$C$}{Set of periods inside a cluster}{}{}. On the basis of these groups, $N_k$ typical or representative periods are derived. 
\begin{figure}[h]
	\centering
  \includegraphics[width=0.49\textwidth]{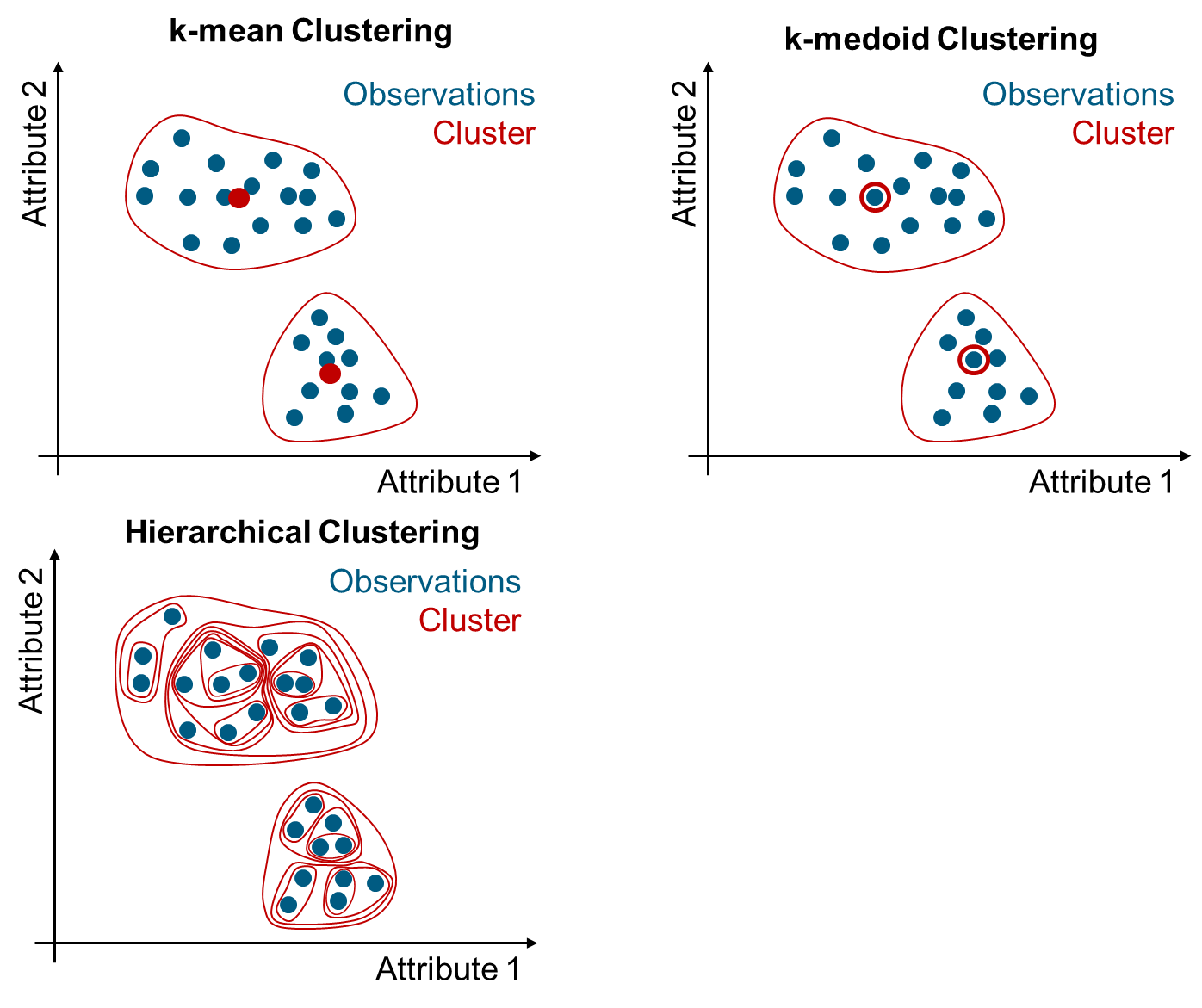}
	\caption{Comparison of cluster methods: k-means with $N_k=2$, k-medoids with $N_k=2$, and hierarchical cluster up to $N_k=2$.}
	\label{fig:ComparisonMeanMedoidHierarchy}
\end{figure}

The representative data point itself of a typical period $k$ \nomenclature[C]{$k$}{Typical period index}{}{}in period step $g$ of attribute $a$ is defined as $\mu_{a,k,g}$\nomenclature[A]{$\mu$}{Representative values of a typical period}{}{}. The representative candidate is stated as vector $\pmb{\mu}_{k}$.
In all reviewed methods, the number of representative periods $N_k$ must be defined a priori and is validated following the aggregation procedure.
The applied cluster methods are visually depicted in Figure \ref{fig:ComparisonMeanMedoidHierarchy}.

\subsubsection{Averaging periods}
Perhaps the most obvious option is to take the averaged values of different periods. \citet{Mavrotas2008} do this for each month in a year for the heating, cooling and electricity load, resulting in 12 typical days.

A more generic formulation that divides the original time series into $N_k$ parts based on their order is introduced as follow:
\begin{enumerate}
\item Calculate the integer divisor of the number of candidate periods to the number of representative periods $\lfloor  \frac{N_i}{N_k}  \rfloor$, which represents the size of each. 
\item Relate $\lfloor  \frac{N_i}{N_k}  \rfloor$ candidate periods in their original order to the cluster period $C_k$, except the last cluster $C_{N_k}$ which receives the remaining candidate periods.
\item Calculate the average period profile of each group as:
\begin{equation}
\mu_{a,g,k} = \frac{1}{\vert C_k \vert}  \sum_{i \in C_k} x_{a,g,i} \quad \forall \quad a,k,g
\end{equation}
\end{enumerate}
The advantage of this method is, that it is easy applicable and the resulting typical periods have a clear order. Nevertheless, the aggregation is only based on the original sequence of the periods and not on the similarity of the group members or candidate periods.

\subsubsection{k-means clustering}

The \textit{k-means} clustering algorithm calculates also the average - the mean - profile of a group. The significant difference to the averaging method is that the k-means algorithm creates the clusters in order to minimize the squared error between the empirical mean of a cluster and all candidates in the cluster \cite{Jain2010}, and not by their original positional appearance in the year. 
Different error or distance functions are possible. In this paper, we consider the squared error, also known as the \textit{Euclidean distance}, which is defined as follows:

\begin{equation}
\min \quad \sum_{k=1}^{N_k} \sum_{i=1}^{N_i} \left[  \sum_{g=1}^{N_g} \sum_{a=1}^{N_a} \left( x_{a,g,i} - \mu_{a,g,k} \right)^2 \right] \times z_{i,k}
\label{eq:objective_clustering}
\end{equation}

or in the vector form as 

\begin{equation}
\min \quad \sum_{k=1}^{N_k} \sum_{i=1}^{N_i}  \left| \left| \pmb{x}_{i} - \pmb{\mu}_{k} \right|\right|^2 \times z_{i,k}
\label{eq:objective_clustering_vec}
\end{equation}

where $z_{i,k}$\nomenclature[A]{$z$}{Binary variable determining candidate cluster assignment}{}{}  is a binary variable that is equal to 1 if the candidate $i$ is assigned to cluster $k$. In order to make sure that each candidate is assigned to a cluster, the following constraint is added:

\begin{equation}
\sum_{k=1}^{N_k} z_{i,k} = 1 \quad \forall \quad i
\end{equation}.

This defines a mixed-integer nonlinear program (MINLP) which is non-convex and difficult to solve. In consequence, the \textit{k-means} method is implemented as a \textit{greedy} algorithm that converges to a local minimum. Nevertheless, it is fairly fast, which is why in practice it is run with many different starting points in order to determine the global minimum.
The main steps of the algorithm itself are as follows:
\begin{enumerate}
\item Randomly or heuristically select an initial partition with $N_k$ clusters and calculate the cluster centers. 
\item Assign each candidate to its closest cluster center.
\item Compute new cluster centers based on all candidates belonging to the cluster.
\item If no convergence criterion is met, return to step 2.
\end{enumerate}

A more detailed discussion of k-means clustering for the selection of typical periods can be found in \citet{Fazlollahi2014}.

\subsubsection{Exact k-medoids clustering}

The k-medoids clustering is an adaption of the k-means algorithm. Instead of calculating the mean as a cluster center, a representative candidate - the medoid - is chosen for $\mathbf{\mu}_k$. 
Figure \ref{fig:ComparisonMeanMedoidHierarchy} portrays this difference. 

Different algorithms for the determination of the clusters also exist. While \citet{Rager2015,Stadler2016} use a greedy optimization algorithm called Partitioning Around Medoids (PAM) \cite{Reynolds2006}, the problem can also be stated as a Mixed-Integer Linear Program (MILP) \cite{Dominguez2011,Schuetz2016}, for which good global solving algorithms exist. 

First, the distance between each candidate is calculated as 
\begin{equation}
d(i,j) =  \sum_{g=1}^{N_g} \sum_{a=1}^{N_a} \left( x_{a,g,i} - x_{a,g,j} \right)^2 \quad \forall \quad i,j \in {1,...,N_i}
\label{eq:distance_candidate}
\end{equation}

Then, following MILP can be stated
\begin{equation}
\min \quad \sum_{i=1}^{N_i} \sum_{j=1}^{N_i}  d(i,j) \times z_{i,j}
\end{equation}
subject to
\begin{equation}
\sum_{j=1}^{N_i} z_{i,j} = 1 \quad \forall \quad j \in {1,...,N_i}
\end{equation}
\begin{equation}
z_{i,j} \leq y_i \quad \forall \quad i,j \in {1,...,N_i}
\end{equation}
\begin{equation}
\sum_{i=1}^{N_i} y_{i} = N_k 
\end{equation}
where both $z_{i,j}$ and $y_i$ \nomenclature[A]{$y$}{Binary variable determining if the candidate is a cluster center}{}{} are binary variables.

This optimization problem can be solved to a global optimum, but has the disadvantage of requiring a high computational load. The computational effort is directly correlated to the number of candidate periods. On the other hand, the impact of the number of attributes and number of time steps in each period is negligibly small, because the distance between the candidates is calculated a priori.

\subsubsection{Hierarchical clustering}

The hierarchical clustering algorithm also minimizes the distance between each candidate period $\pmb{x}_i$ and the representative periods $\pmb{\mu}_k$, as shown in equation \ref{eq:objective_clustering_vec}. The algorithm starts with each single candidate as its own cluster. Then, pairs of clusters are iteratively merged until $N_k$ clusters are left, which is also shown in Figure \ref{fig:ComparisonMeanMedoidHierarchy}. The procedure pf hierarchical clustering is as follows:
\begin{enumerate}
\item Set each candidate as an own cluster.
\item Determine the centroid or mean vector of each cluster.
\item Calculate the \textit{Euclidian} distances to the other cluster centers
\item Merge the two clusters with the lowest distance. If the number of cluster is still greater than $N_k$, return to step 2.
\end{enumerate}
Based on the resulting clusters, representative periods for each of these must be derived. In this work, we therefore chose the \textit{medoid}, respectively the candidate period in the cluster with the smallest distance to all other cluster candidates. In consequence, the method considered is essentially a greedy algorithm of the \textit{k-medoids} algorithm presented in the previous subchapter. 

In comparison to the before presented \textit{k-means} algorithm, the advantage of this algorithm is that it is independent from initial starting points. For this reason, the clusters are easily reproducible.

A more detailed introduction of hierarchical clustering for time series aggregation can be found in \citet{Nahmmacher2016}.



\subsection{Integration of extreme periods}
\label{sec:extreme-periods}

The methods introduced for time series aggregation have the disadvantage of potentially cut off so called peak periods, because such periods are not representative for a whole group or cluster of periods. Rather, they are periods in which the considered time series have design relevant extrema, e.g. peak heat demand. These are important because an accurate energy system design should be able to meet these demands. Therefore, \citet{Stadler2014}, \citet{Dominguez2011}, \citet{Fazlollahi2014}, \citet{Bahl2016} all manually add so called peak periods to the aggregated periods. 

\begin{figure}[h]
	\centering
  \includegraphics[width=0.49\textwidth]{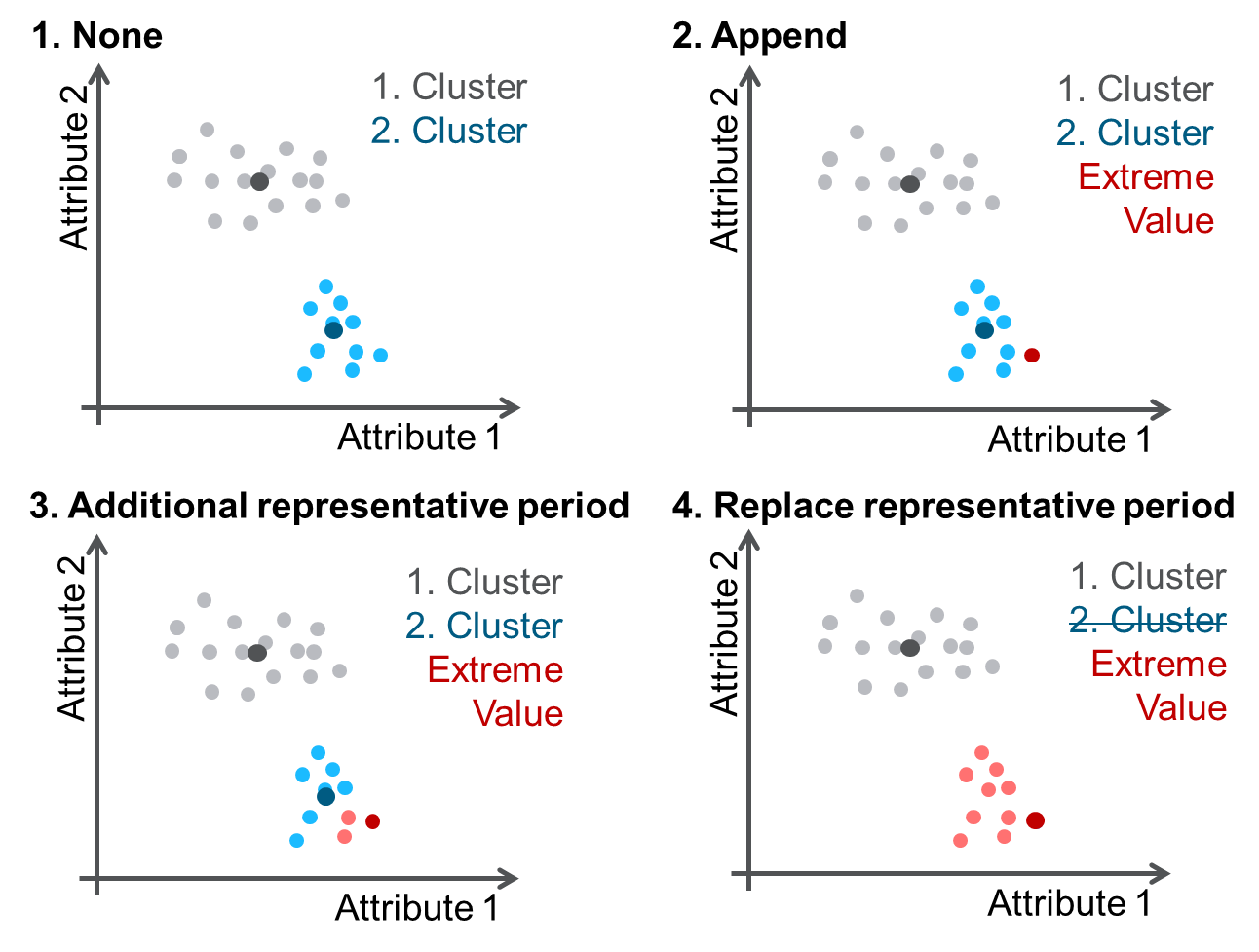}
	\caption{Principle illustration of different methods for integrating the extreme periods, as referred to as peak periods, into the existing aggregated periods.}
	\label{fig:IntegrationExtremePeriods}
\end{figure}

For a comparison of a good integration of these periods, different variants of integration extreme periods are presented in Figure \ref{fig:IntegrationExtremePeriods}. The variants are as follows:
\begin{enumerate}
\item \textit{None} - no integration of extreme periods at all.
\item \textit{Append} - add the extreme periods as additional representative periods to the other representative periods.
\item \textit{Additional cluster center} - the extreme periods are set as an additional new cluster center. Each candidate period is assigned if it is closer to the new cluster center than to its original cluster center.
\item \textit{Replace representative period} - the extreme period is becoming the new representative period of the cluster it was originally assigned to.
\end{enumerate}

The choice of method will determine the robustness of the system design.

\subsection{Scaling of aggregated time-series}
\label{sec:scaling}
In order for the average values of the aggregated time series to fit the average values of the original time series, the time series are scaled a posteriori, similar to the approach of \citet{Dominguez2011,Nahmmacher2016}. It is performed for each time series $a$ separately as follows
\begin{equation}
\mu^{'}_{a,k,g} = \mu_{a,k,g} \frac{\sum_{i=1}^{N_i} \sum_{g=1}^{N_g} x_{a,i,g} }{\sum_{k=1}^{N_k} \vert C_k \vert \times \sum_{g=1}^{N_g}  \mu_{a,k,g} }  \quad \forall \quad a,k,g 
\end{equation}
In order to not exceed the extreme values of the original time series while scaling in accordance to \citet{Nahmmacher2016}, all values greater than 1 are set to 1 and the other values re-scaled again in order to reach the correct average value.

As a last step, the profiles are scaled back to their original scale:
\begin{equation}
\mathbf{\mu}^{''}_{a,k,g} = \mathbf{\mu}^{'}_{a,k,g} \left( \max  \mathbf{x}_{a}^{'} - \min{  \mathbf{x}_{a}^{'} } \right) + \min \mathbf{x}_{a}^{'}   \quad   \forall \quad a,k,g 
\end{equation}

\section{Exemplary time series aggregation}
\label{sec:03_Aggregation}
In validating and comparing the introduced methods, five time series are independently aggregated to four typical days. The time series are hourly values for an entire year. The Global Horizontal solar Irradiation (GHI), the temperature and the wind speed of the test reference year for Bad Marienberg, Germany \cite{DWD2012} outline the weather phenomena relevant to energy system design. Furthermore, two electrical load profiles are aggregated, namely the load profile of a single household \cite{Tjaden2015}, and an electrical load representing an entire region \cite{Robinius2015,Robinius2017}. 

The profiles are reduced to four typical days and analyzed by comparing the original shape and aggregated shape of the profiles (section \ref{sec:qualitativeComparison}). In a further step, indicators are used to contrast the suitability of the different aggregation methods for the different types of profiles (section \ref{sec:IndicaterApplication}).

\subsection{Qualitative comparison of the aggregation methods}
\label{sec:qualitativeComparison}

The number of typical periods is set to $N_k = 4$ and the four aggregation methods described in section \ref{sec:aggregation-methods} are applied.

\begin{figure}[h]
  \includegraphics[width=0.49\textwidth]{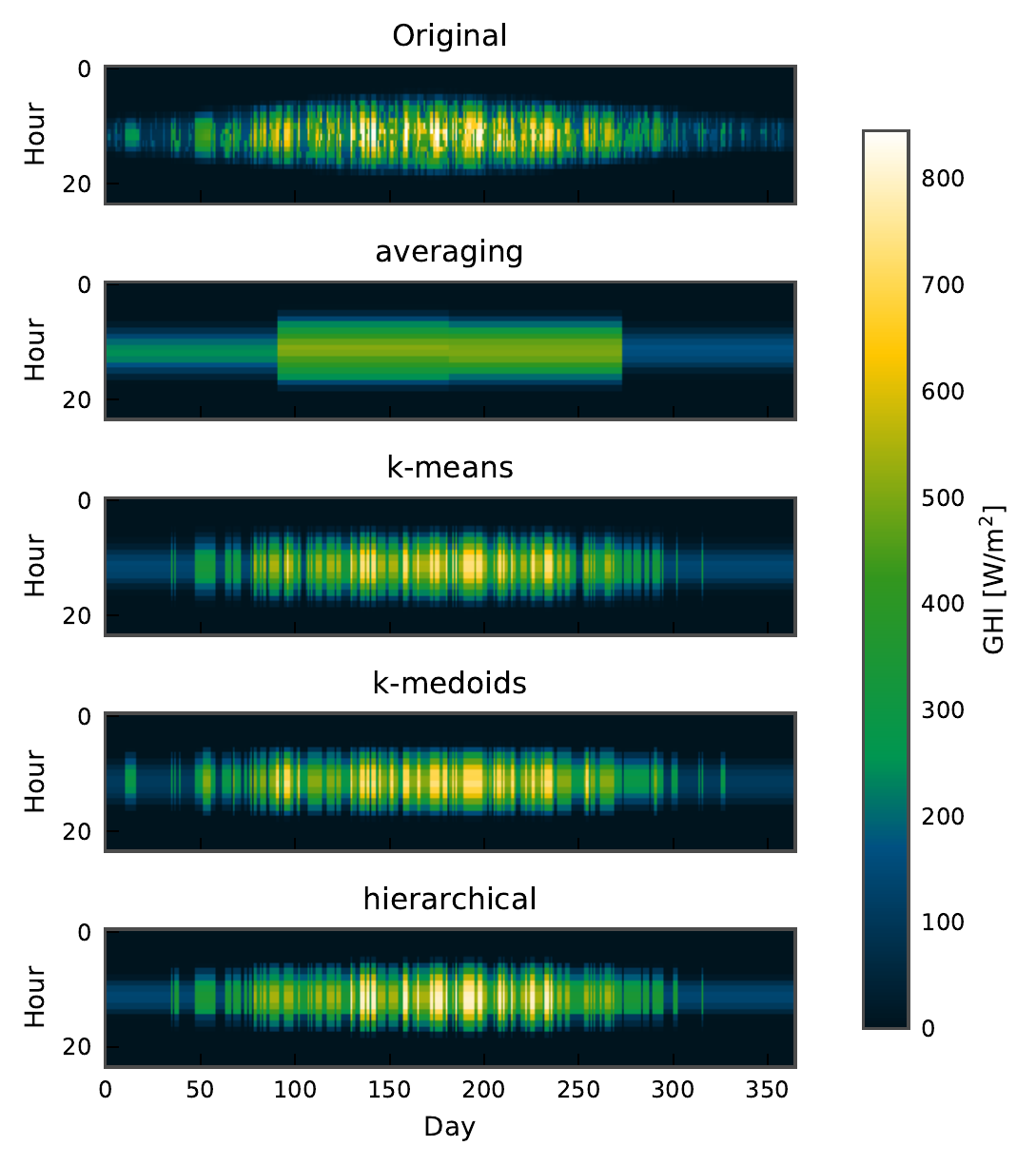}
	\caption{Original annual GHI time series compared to the typical time series with four typical days for different aggregation methods (x-axis: day of the year; y-axis: hour of the day)}
	\label{fig:CMAP_cluster_GHI}
\end{figure}

\begin{figure}[h]
  \includegraphics[width=0.49\textwidth]{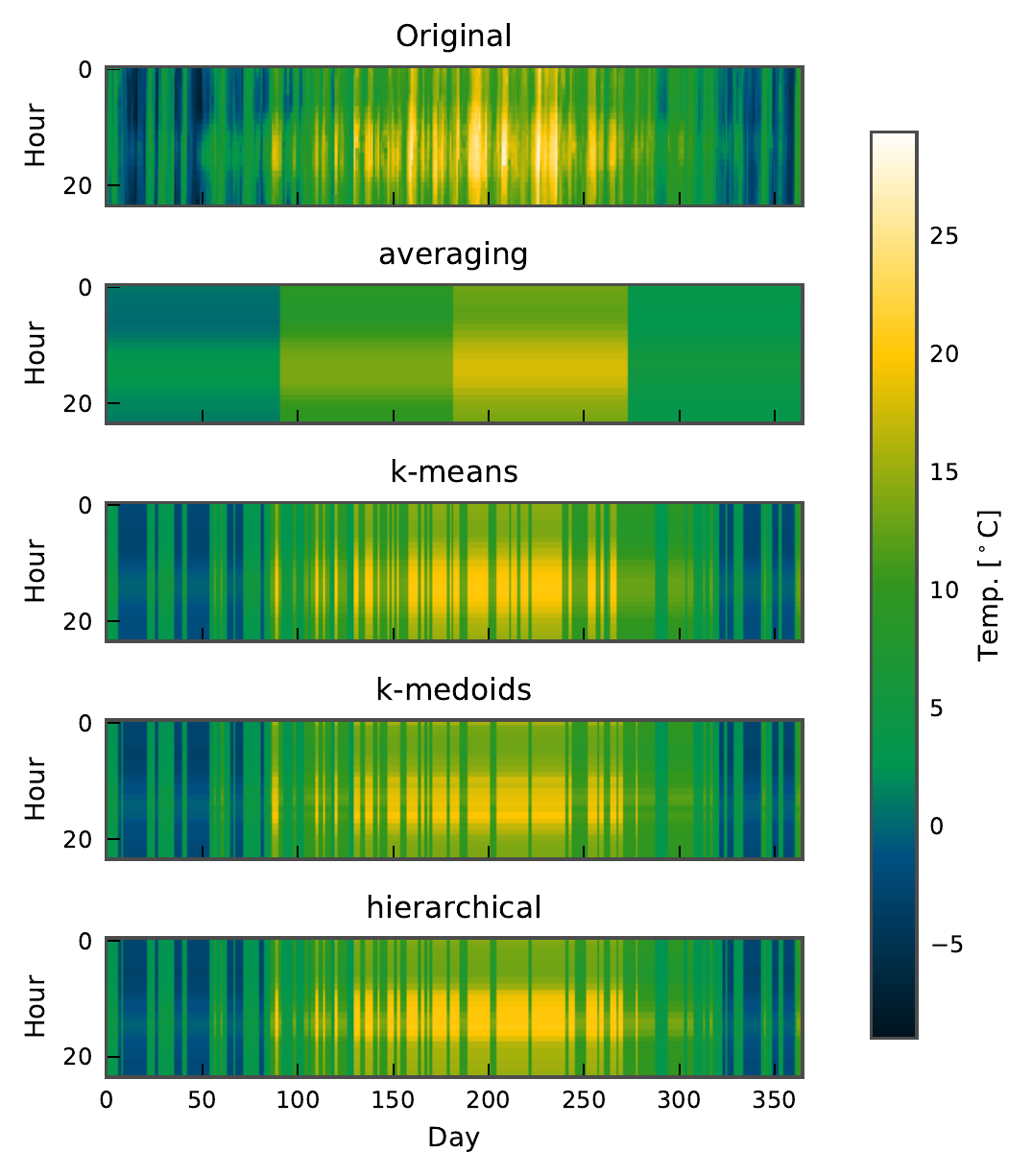}
	\caption{Original annual temperature time series compared to the typical time series with four typical days for different aggregation methods (x-axis: day of the year; y-axis: hour of the day)}
	\label{fig:CMAP_cluster_1}
\end{figure}

\begin{figure}[h]
  \includegraphics[width=0.49\textwidth]{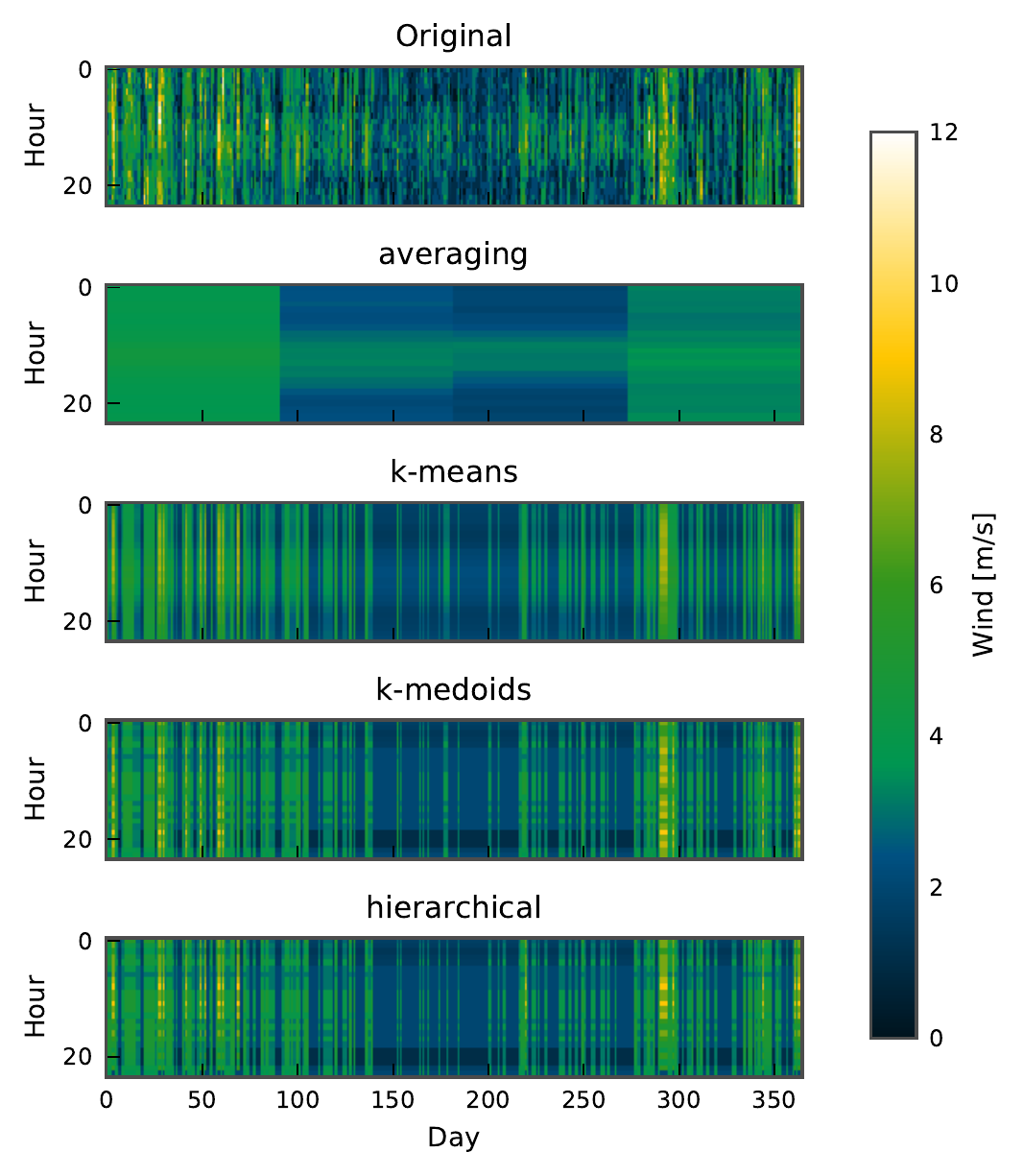}
	\caption{Original annual wind time series compared to the typical time series with four typical days for different aggregation methods (x-axis: Day in year; y-axis: Hour at day)}
	\label{fig:CMAP_cluster_Wind}
  \end{figure}

\begin{figure}[h]
  \includegraphics[width=0.49\textwidth]{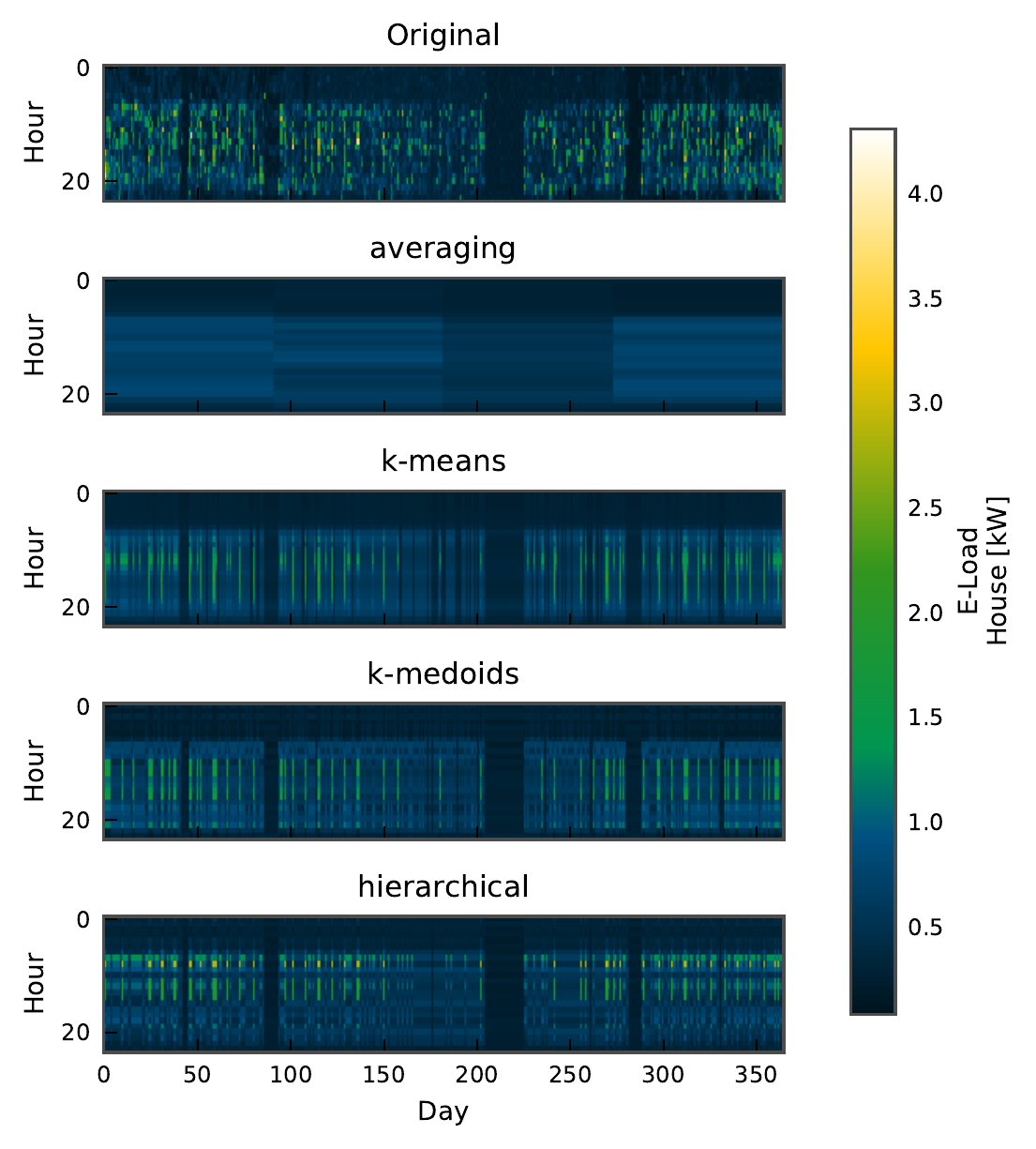}
	\caption{Original annual electrical load time series compared to the typical time series with four typical days for different aggregation methods (x-axis: day of the year; y-axis: hour of the day)}
	\label{fig:CMAP_cluster_2}
  \end{figure}

\begin{figure}[h]
  \includegraphics[width=0.49\textwidth]{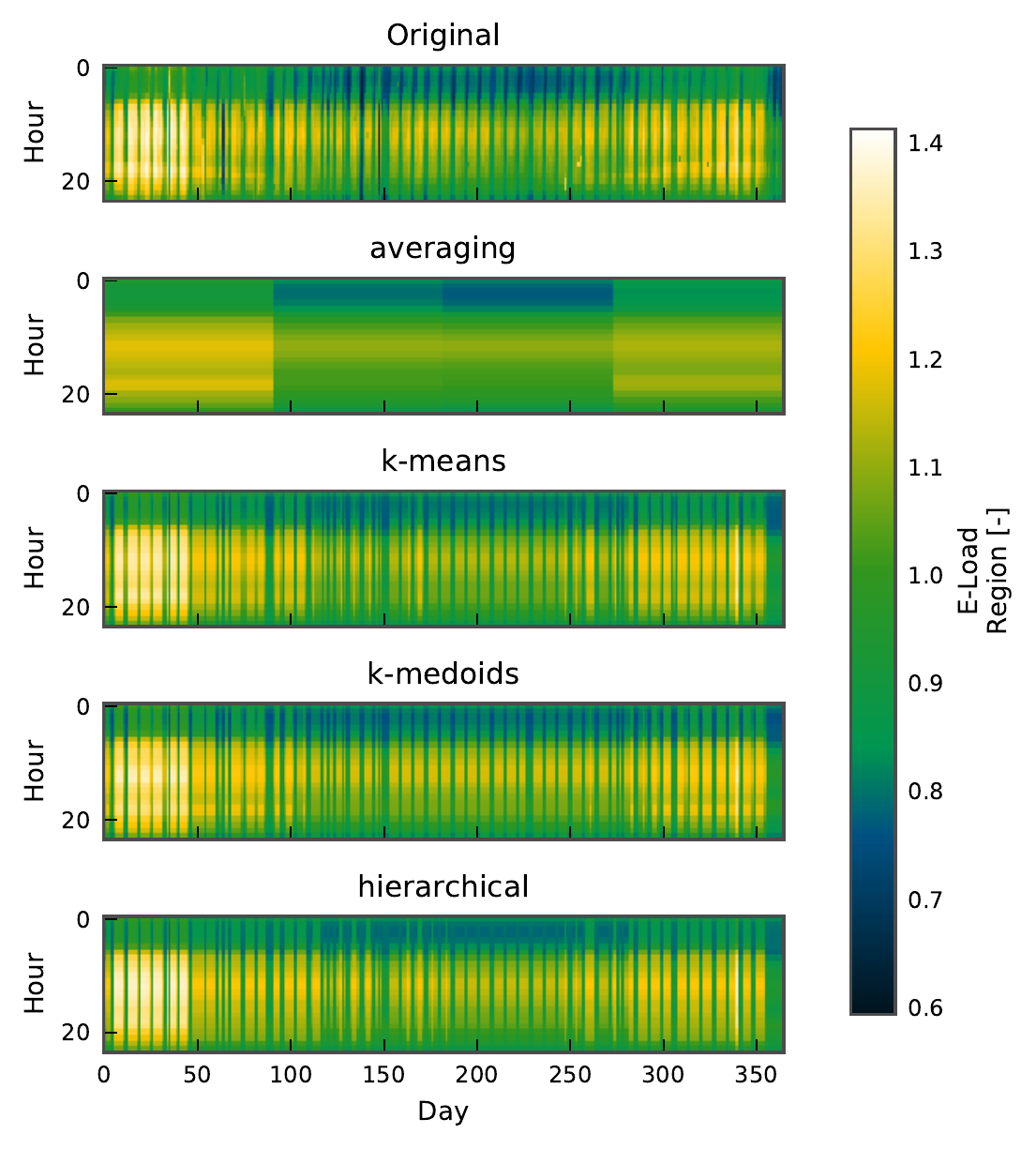}
	\caption{Original regional electricity demand time series compared to the typical time series with four typical days for different aggregation methods (x-axis: day of the year; y-axis: hour of the day)}
	\label{fig:CMAP_cluster_Ele_Region}
\end{figure}

The original profiles are compared to their aggregated counterparts in Figure \ref{fig:CMAP_cluster_GHI} for the GHI, Figure \ref{fig:CMAP_cluster_1} for the temperature, Figure \ref{fig:CMAP_cluster_Wind} for the wind, Figure \ref{fig:CMAP_cluster_2} for the electrical load of a single household and Figure \ref{fig:CMAP_cluster_Ele_Region} for the load of a whole region. 

First, all figures highlight the difference between the aggregation by the clustering algorithms and the aggregation by \textit{averaging} profiles. The averaging method is bound to group a sequence of candidate periods, while the clustering algorithms are able produce aggregated periods which can represent a nonconsecutive order of candidate periods. This degree of freedom allows a more accurate representation of the original time series. The differences between the single groups of the \textit{averaging} method are high, with the consequence that by taking the average of the the candidates, much of the fluctuation in the original periods is smoothed out.

A comparison of the different time series and their aggregation indicates that a high quality discrepancy exists between the different types of profile. 

\textit{Solar radiation} has a strong daily pattern as seen in Figure \ref{fig:CMAP_cluster_GHI} and is fairly well-represented by all clustering-based aggregations. The aggregation of daily profiles benefits from the high daily pattern of the irradiance time series itself, and is able to collect the main daily and annual patterns.

Furthermore, the \textit{temperature} time series in Figure \ref{fig:CMAP_cluster_1} indicates that the major patterns are represented by the aggregated profiles. Nevertheless, a comparison with the original profile indicates that the days with the minimal temperature are not well integrated, although they are highly relevant for energy system design. It highlights why it is important to somehow integrate the extreme periods into the aggregated time series.

\textit{Wind speed }profiles are difficult to aggregate, as seen in Figure \ref{fig:CMAP_cluster_Wind} due to their negligibly small intra-day pattern. The variation between the days is fairly well represented, but the daily profiles themselves seem flattened, as they do not contain the variance of the original profile. Nevertheless, the medoid based methods seem to embody a higher degree of variation than the k-mean clustering.

Figure \ref{fig:CMAP_cluster_2} illustrates the high fluctuation and the lag of a strong pattern of the \textit{electrical load profile for a single household}. In consequence, the aggregation methods are unable to represent it well and cut out many of the peak loads. Still, a quality difference is once again recognizable between the mean based aggregation and the medoid based aggregation, because the chosen centroid of the k-mean is smoothing out more of the fluctuations. Nevertheless, for this profile the results also indicate that for all chosen methods, it is necessary to add the related peak periods afterwards.

In contrast, Figure \ref{fig:CMAP_cluster_Ele_Region} shows the \textit{electrical load of an entire region} and its strong daily pattern. Further differences between winter and summer, as well as between weekdays and weekends, are apparent. Due to these strong patterns, all three cluster methods can aggregate the design-relevant daily profiles, while no qualitative difference is recognizable between the different aggregations.

For all cases, it is noticeable that the assignment of candidate periods to typical periods - the grouping - does not variate much between the three different clustering algorithms. 

\subsection{Indicator based quality evaluation}
\label{sec:IndicaterApplication}
In the literature many different indicators have been introduced to evaluate which aggregation is the most suitable \cite{Poncelet2016,Fazlollahi2014}. The use of indicators has the advantage that they allow an evaluation of the accuracy of the time series aggregation before any optimization procedure has even started. 
Two indicators have been used in this work: 
\begin{itemize}
\item The Root-Mean-Squared-Error (RMSE)\nomenclature[D]{RMSE}{Root Mean Squared Error}{}{} which is equivalent to the average intra-cluster distance introduced in equation \ref{eq:objective_clustering}. Essentially, it is the objective functional value of all aggregation methods.
\item The RMSE between the duration curve - the sorted data set - of the original time series and the aggregated time series. The comparison of the duration curves between the aggregated and original profile is important because it allows the evaluation of all potential data points, which are represented independently. 
\end{itemize}

The first indicator (RMSE) for the normalized profiles for the different profiles and for the different aggregation methods is shown in Table \ref{tab:RMSE} which displays the same data as was presented in the previous section.

\begin{table}[h]
	\centering
\caption{Root-Mean-Squared-Error (RMSE) in percent between the original profile and aggregated profile of four typical days for different aggregation methods and different types of time-series}
\resizebox{\columnwidth}{!}{%
\label{tab:RMSE}
\begin{tabular}{lccccc} 
\\ \hline 
 &  GHI  &  Temp. &  Wind &  E-Load House &  E-Load Region \\ \hline averaging    &          10.45 &              11.57 &       14.44 &                9.39 &               12.90 \\ k-means      &           6.45 &               6.32 &        9.99 &                7.75 &                6.02 \\ k-medoids    &           6.83 &               6.55 &       10.78 &                8.37 &                6.32 \\ hierarchical &           7.15 &               6.59 &       11.00 &                9.25 &                6.55 \\ 
\end{tabular}
}
\end{table}

The comparison of the RMSE between the profile types supports the interpretation of the qualitative comparison concerning the quality of the  aggregation: The electricity load of a whole region has the smallest RMSE, at 6.02 \%  average for the k-mean clustering, followed by temperature and GHI. With a big gap, the electrical load of a household follows with an error of 7.75 \%, while the wind profile has the highest value, with an error of 9.99 \%.

The comparison in terms of the aggregation methods shows the following:
The \textit{averaging} method performs the worst for all profiles, the \textit{k-mean} clustering performs best at the  RMSE in comparison to the \textit{k-medoids} and \textit{hierarchical} clustering. This is because a centroid generally features less distance to its candidates in comparison to a medoid. The \textit{hierarchical} aggregation performs the worst in all cases. 

Nevertheless, a look at the RMSE of the duration curves in Table \ref{tab:RMSE_duration} shows a less uniformly distribution between the profile types and the aggregation methods. The \textit{hierarchical} aggregation performs the best for the aggregation of GHI with 1.71 \% mean error, and with the electrical load of a single house featuring a 3.06 \% error. \textit{k-medoids} get the best results for the wind speed with 4.09 \% error and the electrical load of a region with 1.76 \% error. The temperature is best represented in terms of the duration curve by the \textit{k-means} aggregation, with an error of 2.65 \%.

\begin{table}[h]
	\centering
\caption{Root-Mean-Squared-Error (RMSE) between the duration curve of the original profile and that of the aggregated profile of four typical days for different aggregation methods and different profile types}
\label{tab:RMSE_duration}
\resizebox{\columnwidth}{!}{%
\begin{tabular}{lccccc} 
\\ \hline 
 &  GHI  &  Temp. &  Wind &  E-Load House &  E-Load Region \\ \hline averaging    &           5.10 &               5.55 &       10.43 &                7.63 &                6.11 \\k-means      &           1.93 &               2.65 &        4.63 &                3.95 &                1.94 \\ k-medoids    &           1.75 &               3.21 &        4.09 &                4.05 &                1.76 \\ hierarchical &           1.71 &               2.74 &        4.28 &                3.06 &                2.57 \\ 
\end{tabular}
}
\end{table}

The hierarchical and k-medoid aggregation perform better in case of the RMSE of the duration curves because the medoids represent more fluctations and also more peaks in comparison to the centroids, which partially average them  out. These peaks are also represented in the original duration curve, which is why these aggregations have a higher overlap.

Still, it would be difficult to favor a certain method based on the indicators introduced. Their information value is limited and it is hard to predict how they would affect optimal system design. To overcome this, we will apply the typical periods in the following sections to the design optimization of different energy systems .

\section{Application of optimal energy system design}
\label{sec:examples}
For analyzing the advantages and drawbacks of time series aggregation for energy system optimization, the methods introduced are applied to three different energy supply systems:
\begin{enumerate}
\item A combined heat and power plant (\textit{CHP})\nomenclature[D]{CHP}{Combined Heat and Power plant}{}{} in combination with a gas boiler and heat storage for the supply of the residential electricity and heat demand, introduced in section \ref{sec:CHPsys}.
\item A \textit{residential} supply system that is largely based on a heatpump and photovoltaics, introduced in section \ref{sec:ResSystem}. 
\item An \textit{island} system that supplies the electricity demand of an entire region, as introduced in section \ref{sec:IslandSystem}. This is based on photovoltaics and wind as power sources and also needs storage technologies because it has a limited access to fossil resources. 
\end{enumerate}

The aggregation methods and number of typical periods are varied and compared for all three systems. Furthermore, in the case of the residential system, we compare in detail the impact of different methods on integrating extreme periods. The island system is introduced to compare the impact of the varied length of the typical periods - e.g. typical days or typical weeks - on the system's design. 

The operation and design of the systems are optimized to achieve minimal costs of the energy supply. The general modeling approach of the systems is introduced in \ref{sec:systemModel}. The models are fairly simple in order to achieve a rapid solving performance in producing results for many different aggregation variants. 
The modeling language is \textit{Pyomo 4.3}  \cite{Hart2011} and as solver \textit{Gurobi 7.0.1}  \cite{Gurobi2016} was chosen. The hardware was an Intel i7-4790 CPU with 32 GB RAM, where 7 of 8 threads were used for solving.

\subsection{CHP-based supply system}
\label{sec:CHPsys}

\begin{figure}[h]
	\centering
  \includegraphics[width=\columnwidth]{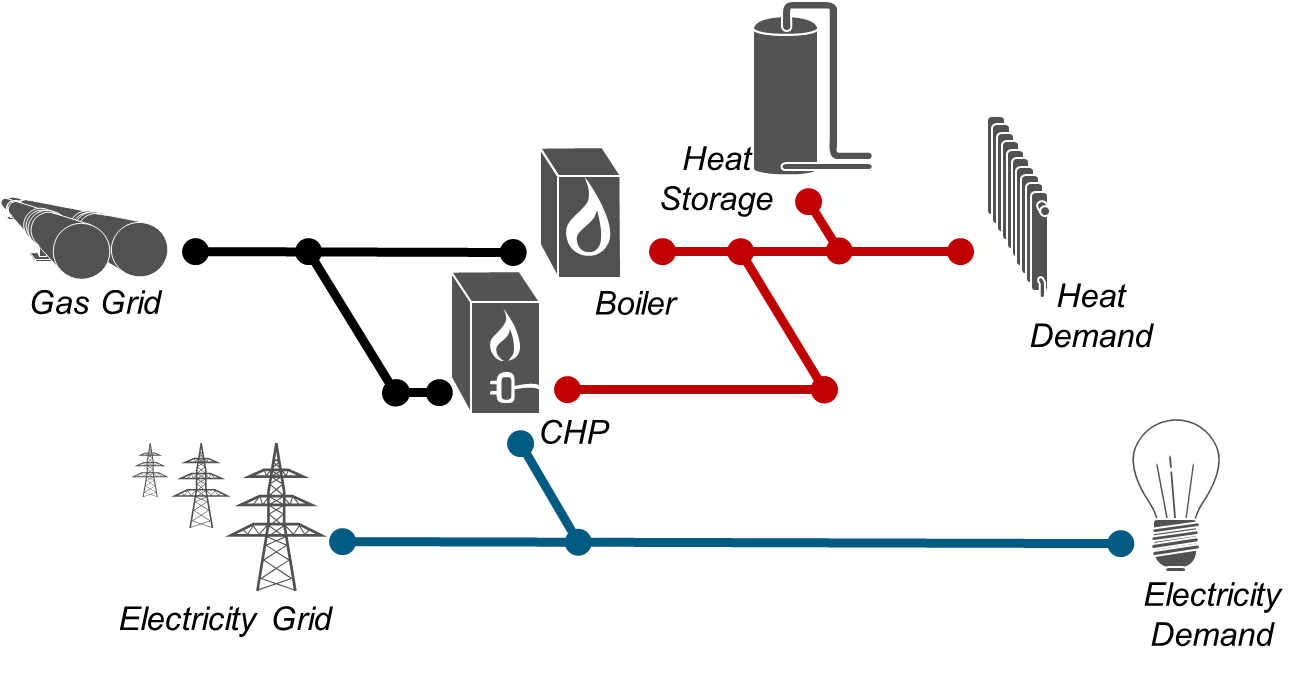}
	\caption{Component network of the considered CHP based energy supply system}
	\label{fig:CHP_System}
\end{figure}

Optimizing the configuration and operation of a CHP system is a common application for typical period aggregation, especially if integer variables are included for choosing real components or modeling discrete states in the operation of the system. 
The CHP system is shown in Figure \ref{fig:CHP_System} and adds to the CHP itself a peak boiler and heat storage, which are all getting scaled and operated. The heat demand is simulated by a 5R1C model \cite{Schuetz2017,EN2008}, with the introduced weather data as an input for a multi-family house. The electrical load is aggregated by six single household profiles that \citet{Tjaden2015} represent the housing units in the multi-family house.

\subsubsection{Comparing the aggregation methods}

For the validation, we optimized the system for the full original time series, and subsequently for different numbers of typical days, aggregated by the methods shown in section \ref{sec:aggregation-methods}. The integration of peak periods was achieved by setting the days with the highest heat and highest electricity demand as \textit{additional cluster centers}, which was discussed in section \ref{sec:extreme-periods}.

\begin{figure}[h]
	\centering
  \includegraphics[width=\columnwidth]{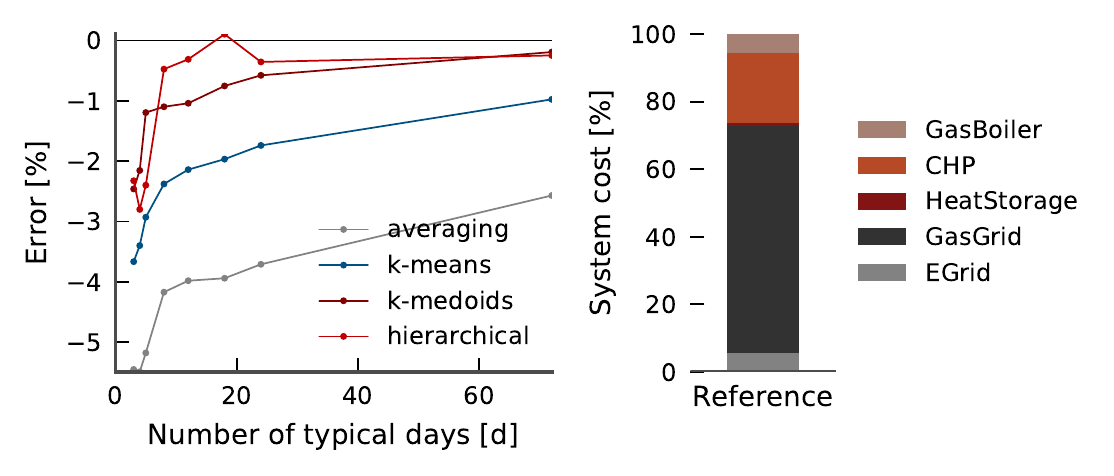}
	\caption{Relative error of the objective value, the annual energy cost, of the optimal CHP-system design based on aggregated periods for different types of aggregation methods in comparison to the annual cost of the cost optimal CHP-system design based on the full time series as reference.}
	\label{fig:Error_Min_CHP}
\end{figure}

The results can be seen in Figure \ref{fig:Error_Min_CHP}. It shows the relative error between the objective value, the annual cost, of the cost-optimal CHP system design based on aggregated periods, in comparison to the annual cost of the cost-optimal CHP system design based on the full time series for different types of aggregation methods. The bar plot is the composition of the annual cost of the optimal system design based on the full time series.

It is noticeable, that all aggregation methods underestimate the cost for a small number of typical periods, although the extreme periods are included. The reason for this is in the approach itself: The most representative periods are in general the days which have a smaller variance and smoother profiles than the original time series. The consequence is that the typical heat demand and the typical electricity demand fit better to each other, and the demand values do not exceed that often the capacity of the CHP plant. This results in an increased supply rates by the CHP system, with the consequence that for example the coverage with the CHP-plant is overestimated and the resulting system cost underestimated.

Therefore, it is also reasonable that the medoid based methods (\textit{k-medoids} and \textit{hierarchical}) perform better than the centroid based method \textit{k-mean}. The choice of real days includes higher variations in comparison to an averaged profile as in the \textit{k-mean} method. The \textit{averaging} method performs, as expected, the worst because it is smoothing out to many fluctuations. No dominance of the \textit{hierarchical} clustering or either the \textit{k-medoid} clustering exists. For both methods it is possible to reach errors in the annual cost smaller than 2\% with 8 typical days. 

\subsubsection{Gain in computational solving duration}
\begin{figure}[h]
	\centering
  \includegraphics[width=\columnwidth]{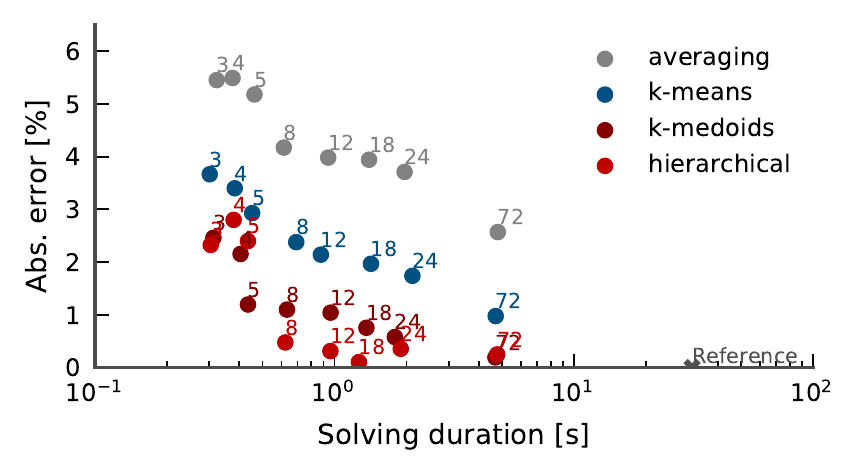}
	\caption{Trade-off between the error in the objective function and solving duration for different numbers of typical days, and aggregation methods}
	\label{fig:ErrorVSPerformance_Min_CHP}
\end{figure}

The gain in the solving performance as a trade-off to the related accepted error is illustrated in Figure \ref{fig:ErrorVSPerformance_Min_CHP}. As expected, the reduction in time steps results in a reduction of the solving time. The shape relates to the type of system model, while for our simple CHP-MILP we can consider a roughly logarithmic reduction of the error related to the solving performance. The precise shape depends on the scaling of the optimization problem and is not generalizable. Nevertheless, this pareto-front can be drawn for each system and it is up to the user to decide which degree of accuracy should be chosen. 
For the case of eight typical days, the solving time could get approximately reduced by a factor of 50 in solving based on the annual profile.

\subsection{Residential supply system}
\label{sec:ResSystem}
\begin{figure}[h]
	\centering
  \includegraphics[width=\columnwidth]{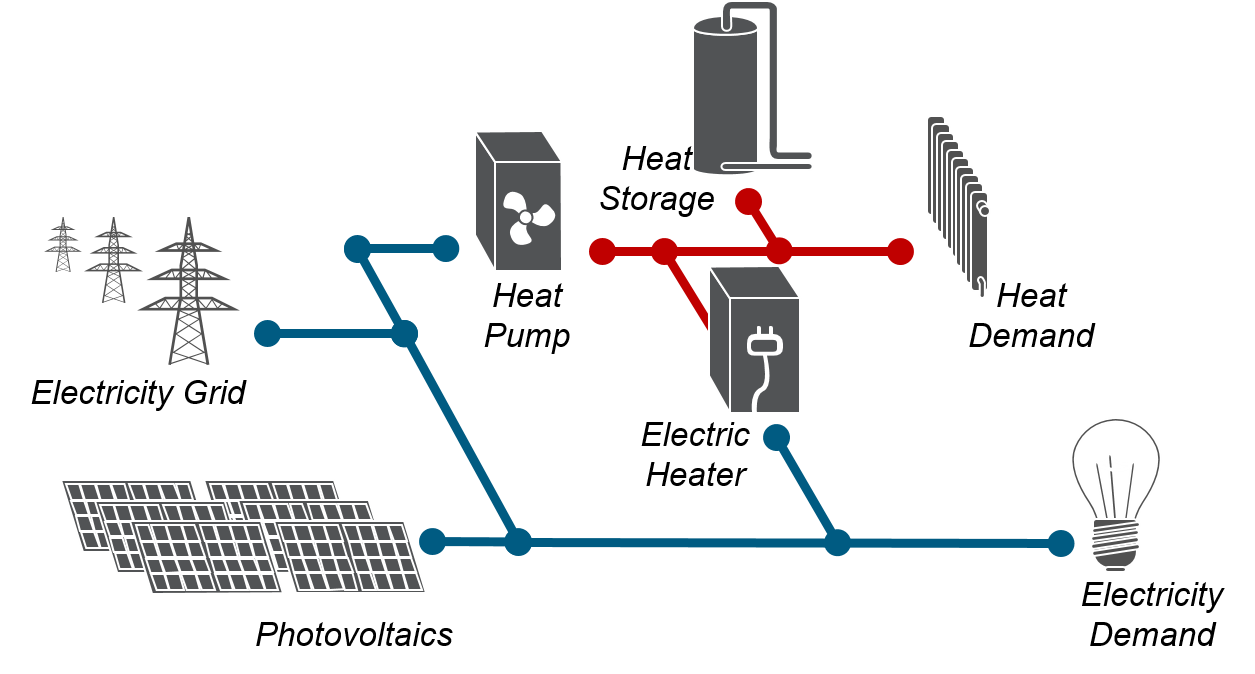}
	\caption{Component network of the residential energy supply system considered}
	\label{fig:PV_HP_System}
\end{figure}

The second system analyzed is a supply system for a residential building based on a central grid supply and photovoltaics. The heat is supplied by a heat pump, an immersion heater with heat storage, as is seen in Figure \ref{fig:PV_HP_System}. Heat and electricity loads are given for a single family house, while the photovoltaic feed-in profile is simulated in advance with the PV-Lib \cite{Andrews2014}.


\subsubsection{Comparing the aggregation methods}
\begin{figure}[h]
	\centering
  \includegraphics[width=\columnwidth]{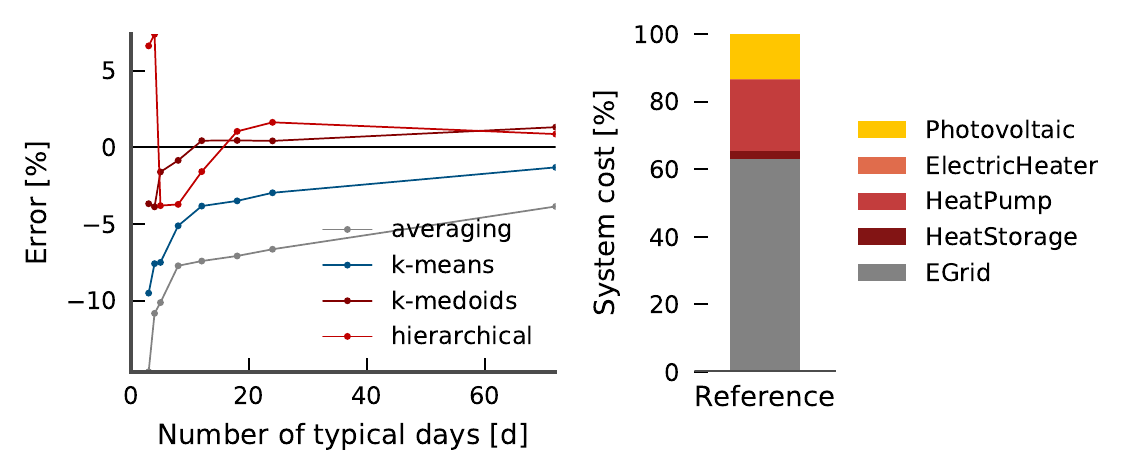}
	\caption{Relative error of the objective value, the annual energy cost, of the optimal residential system design based on aggregated periods in comparison to the annual cost of the cost optimal residential system design based on the full time series for different types of aggregation methods.}
	\label{fig:Error_Min_HP_PV_Cost}
\end{figure}

Like the CHP-system, different aggregation methods have been used for creating the typical periods, while varying also the number of typical days. The results are shown in Figure \ref{fig:Error_Min_HP_PV_Cost}.

It is noticeable that the shape of the error functions of the residential systems are similar to the shape of the CHP system's error function. Furthermore, for the residential system, \textit{k-mean} clustering and \textit{averaging} perform poorly in terms of underestimating the system cost, while the medoid based aggregation methods converge quickly on costs which are similar to the cost of the reference system, but overestimate the cost for higher numbers of typical days.

The high overestimates of the system cost of the \textit{hierarchical} aggregation of three and four typical days show the sensitivity to the actual choice of representative days, e.g. a medoid day has a small overlay of the photovoltaic feed-in and electrical load profiles, self-supply rates are underestimated, with the effect that the energy cost are overestimated. These effects appear less with the \textit{k-mean} aggregation because its profiles are averaged. Therefore, they converge more smoothly on smaller errors with an increased number of days.

The scale of the errors is in average twice as high as the error of the CHP systems for the same number of typical days. Therefore two factors must be taken into account:
\begin{enumerate}
\item The electricity demand of a single household shows a higher fluctuation than that by an aggregated number of households, which is more difficult to capture in typical periods.
\item An additional profile - the photovoltaic feed-in - impacts the optimal system design, which increases the number of patterns to collect with the aggregation methods.
\end{enumerate}
Nevertheless, only 12 typical days are enough for reaching an error, corresponding to less than 2\% of the objective value.


\subsubsection{Comparing the integration of peak periods}
In order to find the best method for the integration of peak periods, we applied all the methods introduced in section \ref{sec:extreme-periods} to the residential system designed by six typical days aggregated by the \textit{k-medoid} clustering.
We manually defined for which time series, what type of peak period should be integrated:
\begin{itemize}
\item the day with peak heat load
\item the day with peak electricity load
\item the day with the smallest total photovoltaic feed-in
\end{itemize}

\begin{figure}[h]
	\centering
  \includegraphics[width=\columnwidth]{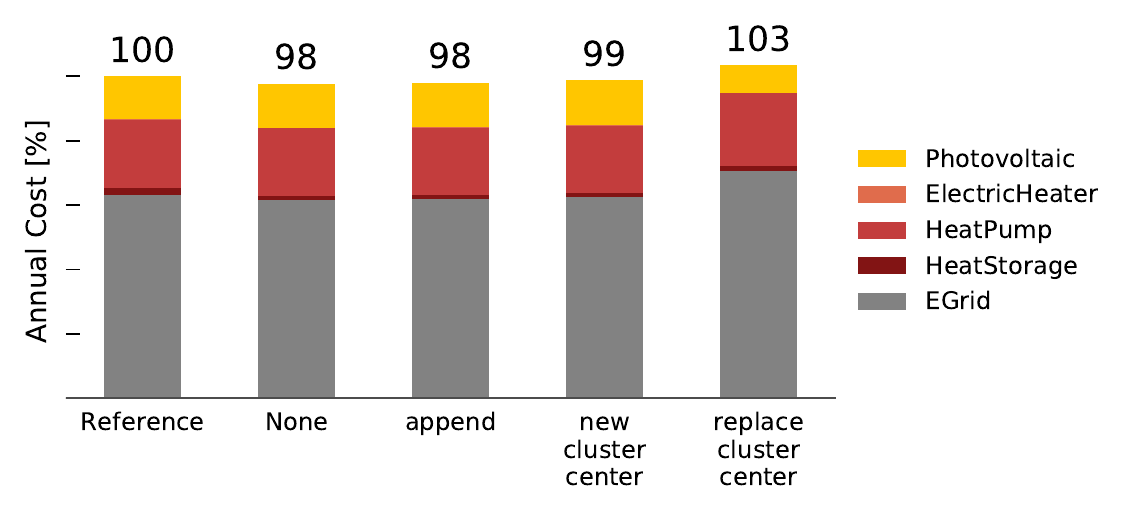}
	\caption{Composition of the system cost for an optimal residential system design for the full time series and the optimal system design based on six typical days with differing integrations of extreme days}
	\label{fig:Error_IntegrationMethod}
\end{figure}

\begin{figure}[h]
	\centering
  \includegraphics[width=\columnwidth]{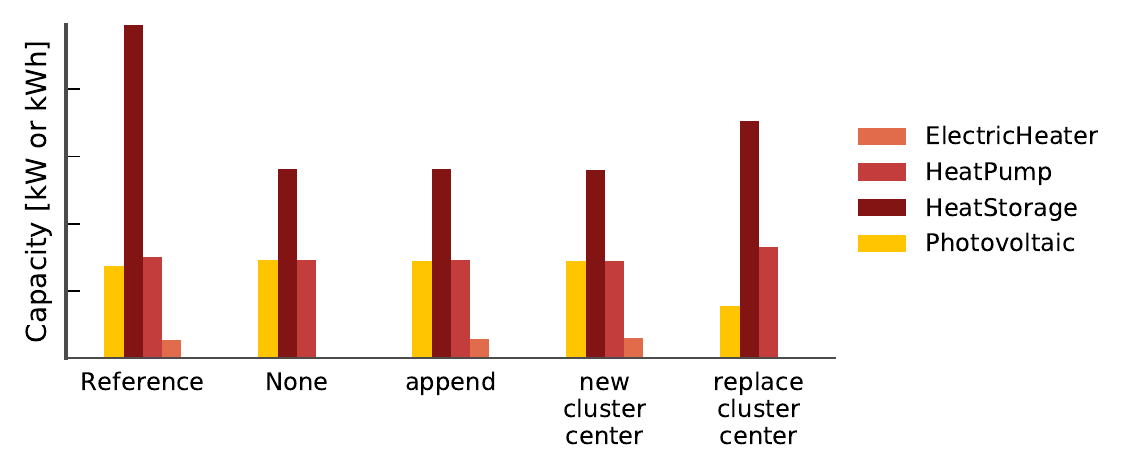}
	\caption{Scaling of the system components for an optimal residential system design for the full time series and system design based on six typical days with different integrations of extreme days}
	\label{fig:integrating_extreme_periods_capacity}
\end{figure}

The results can be seen in Figure \ref{fig:Error_IntegrationMethod}. If no extreme periods are integrated, the system costs are underestimated by roughly 2\%. No improvement is recognizable apart from appending the peak periods. Setting the peak periods manually as a potential cluster center, the underestimation drops to 1\%. The resulting design difference between these two integration methods is so small due to the fact that the peak periods are mostly at the edge of a cluster. In consequence, no or not many other observational periods are related to them. 
The last method, of replacing a whole cluster by a peak period, results in an overestimation of 3\% and forces a more conservative system design because extreme periods are represented above average. A consequence of this is that the system design relies more on electricity supplied by the grid and less on photovoltaics as seen in Figure \ref{fig:integrating_extreme_periods_capacity}.

Also noteworthy is that for all the system designs derived from typical periods, the optimal photovoltaic capacity is overestimated, while the optimal storage capacities get underestimated. This is because of the already mentioned smoothing effect in the typical periods. The photovoltaic feed-in and electricity demand are predicted to be more steady, whereby a higher overlap exists. Furthermore, the heat demand is flattened, with less heat storage capacity required to balance it.


\subsection{Island system}
\label{sec:IslandSystem}

\begin{figure}[h]
	\centering
  \includegraphics[width=\columnwidth]{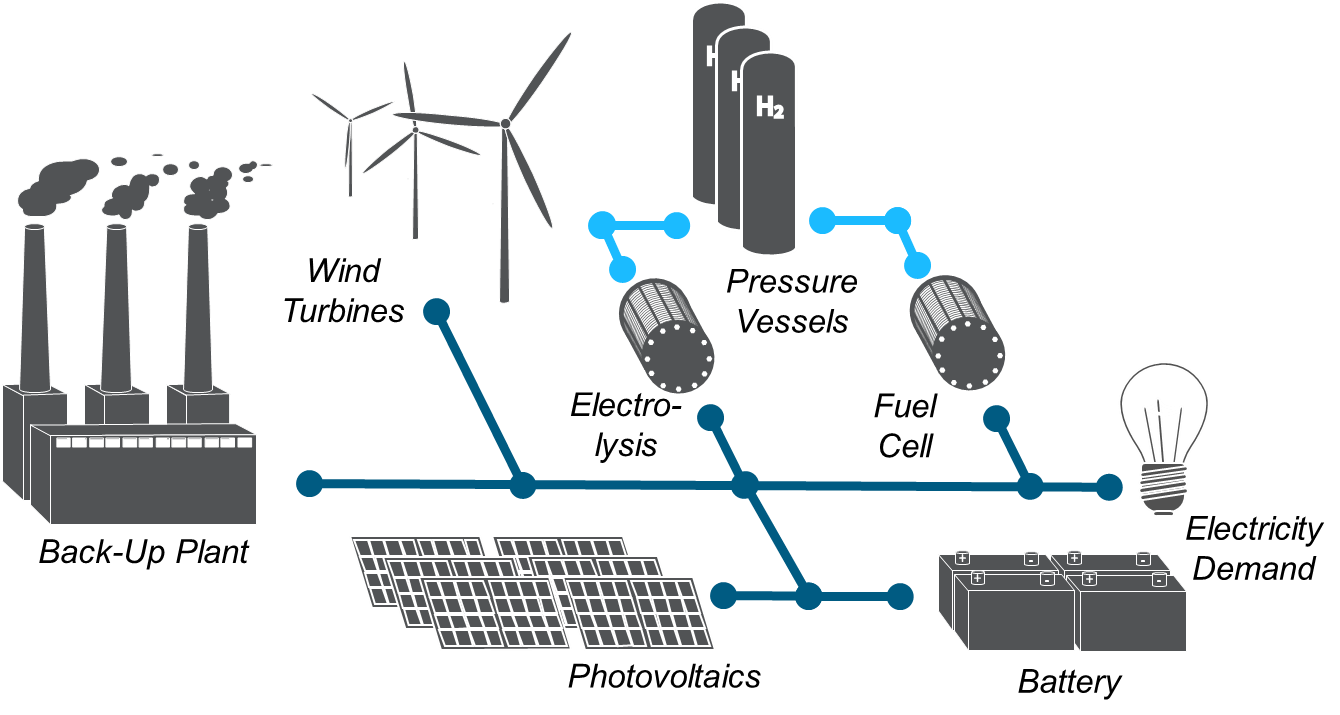}
	\caption{Component network of the island electricity supply system considered}
	\label{fig:Wind_System}
\end{figure}

In terms of validating time series aggregation to a highly renewably based system design, we introduce an island system that is largely reliant on wind turbines and photovoltaics. The regional electricity demand profile and wind feed-in are drawn from \citet{Robinius2017}. The potential photovoltaic feed-in is also simulated with PV-Lib \cite{Andrews2014}. A power plant gets scaled as backup capacity while its total feed-in is restricted to 10\% of the total electricity demand, so as to achieve a high share of renewable energy supply. Furthermore, hydrogen storage based on electrolysers \cite{Schiebahn2015}, fuel cells and pressure vessels, along with a battery storage can balance the fluctuating feed-in of the renewable time series. The entire system is shown in Figure \ref{fig:Wind_System}.

\subsubsection{Comparing the aggregation methods}
\begin{figure}[h]
	\centering
  \includegraphics[width=\columnwidth]{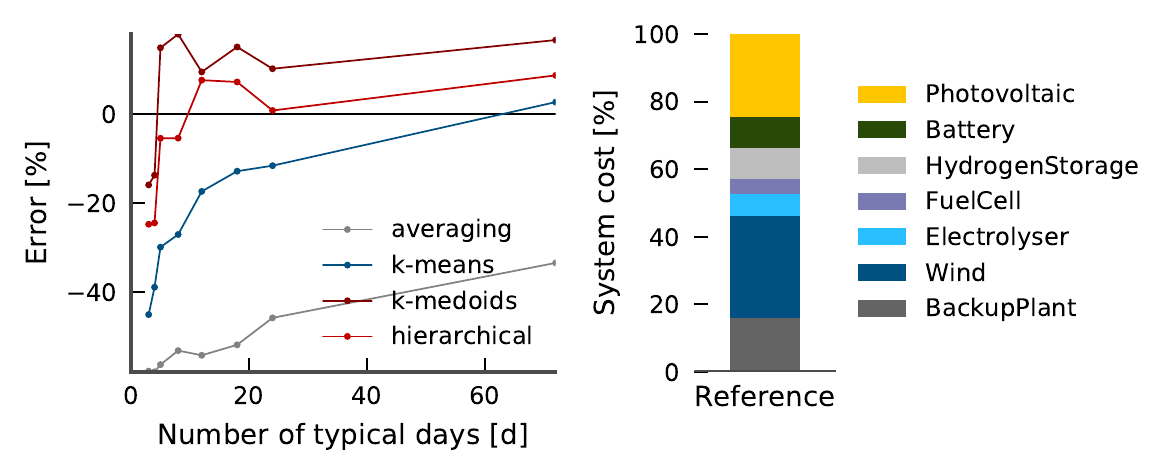}
	\caption{Relative error of the objective value, the annual energy cost, of the optimal island system design based on aggregated periods in comparison to the annual cost of the cost optimal island system design based on the full time series for different types of aggregation methods.}
	\label{fig:Error_Min_Wind_Cost}
\end{figure}

As with both prior systems, the different aggregation methods and the number of typical days are first varied. The resulting optimal system design and errors are depicted in Figure \ref{fig:Error_Min_Wind_Cost}.

The results in terms of the aggregation methods are the same as for both systems introduced previously. The major difference is the scale of the error: While the CHP system and the residential system both converged within less than 12 typical days to an error of less than 2\% for the medoid-based aggregation methods, even for 72 typical days no robust convergence is identifiable for the island system. 

As in the afore-mentioned systems, the smoothing-effect of the aggregation also partially contributes to the error, but therefore an underestimation of the resulting system cost would be expected. However, \textit{k-medoid} and \textit{hierarchical} clustering overestimate the system's costs especially for high numbers of typical days.

Therefore, we expect the cause of the error to lie elsewhere: The cost optimal system design relies heavily on storage systems that contribute to almost 30\% of the system costs, as is seen in Figure \ref{fig:Error_Min_Wind_Cost}, including fuel cell, electrolyzer, hydrogen storage and battery. With the applied approach, each typical period is considered independently and no energy exchange between the typical periods is possible. For our case, this means that only storage operation within a single day is possible.

\subsubsection{Variation of the period length}
For this reason, we vary the duration of the typical periods. Along with the typical days ($N_k=24$), typical groups of three days ($N_k=72$) and typical weeks ($N_k=168$) are also aggregated. We afix the aggregation methods to \textit{k-medoids} and integrate the extreme periods as an additional cluster center. 

\begin{figure}[h]
	\centering
  \includegraphics[width=\columnwidth]{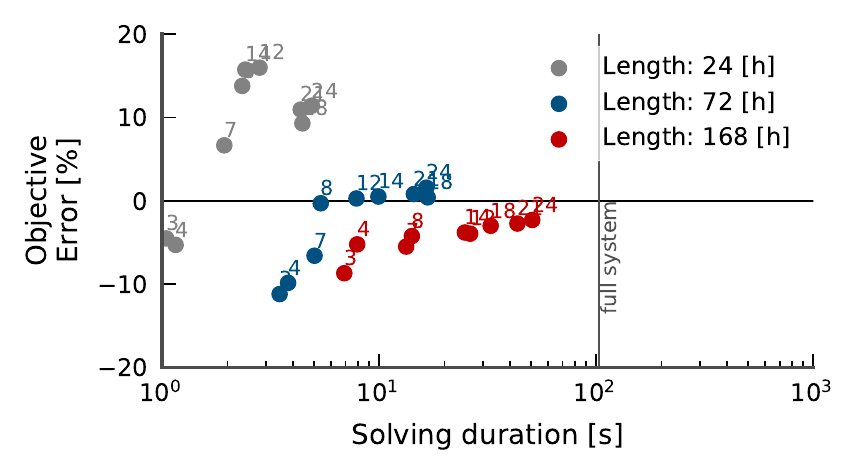}
	\caption{Trade-off between the error in the objective function and solving duration for different numbers of typical periods, and for different period lengths}
	\label{fig:ErrorVSPerformance_PeriodLength_Wind}
\end{figure}

Figure \ref{fig:ErrorVSPerformance_PeriodLength_Wind} shows the resulting error in the objective function and solving duration for the three different period lengths and different numbers of typical periods. 

As was seen before, for typical days almost no convergence is recognizable. 

The aggregation by groups of three days performs slightly better: With higher numbers of typical periods, the error gets reduced, but this convergence also appears non-monotonuously. 

More predictable results can be achieved with typical weeks, which converge continuously on smaller errors. Nevertheless, the error also remains over 2\% for 24 typical weeks. Furthermore, their solving performance is also worse due to the higher number of time steps per period.

\begin{figure}[h]
	\centering
  \includegraphics[width=\columnwidth]{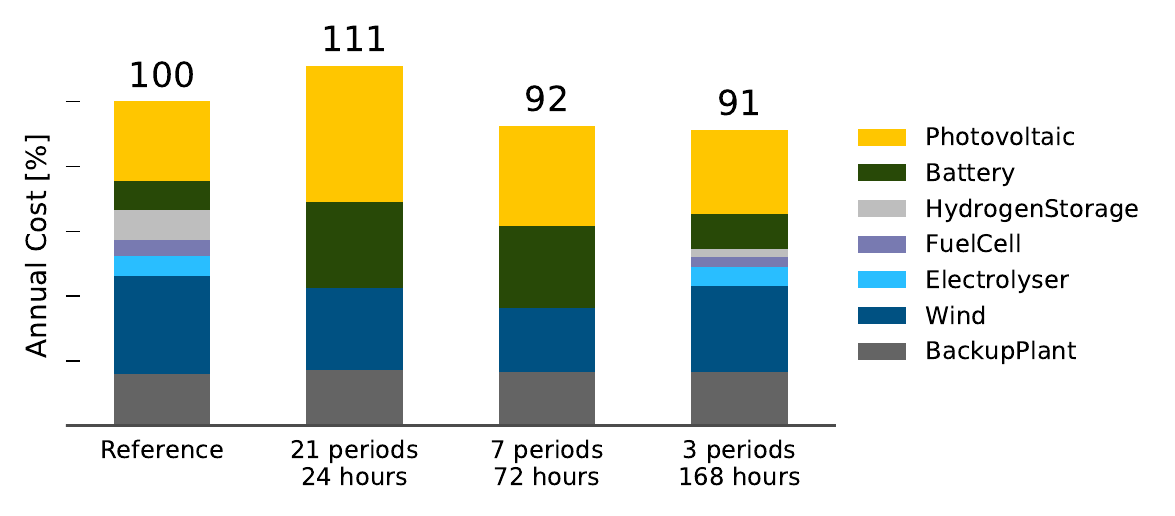}
	\caption{Composition of the system cost for an optimal island system design for the full time series and the optimal system design based on 21 typical days, seven groups of three days, and three typical weeks}
	\label{fig:Error_IntegrationMethod_Periods}
\end{figure}

For comparing aggregations with comparable solving performance, Figure \ref{fig:Error_IntegrationMethod_Periods} shows the cost composition for the aggregation of 21 typical days (each 24 time steps), seven typical periods based on three days (each 72 time steps) and three typical weeks (each 168 steps), which in sum have all the same number of time steps. They are compared to the original system design for validation.

\begin{figure}[h]
	\centering
  \includegraphics[width=\columnwidth]{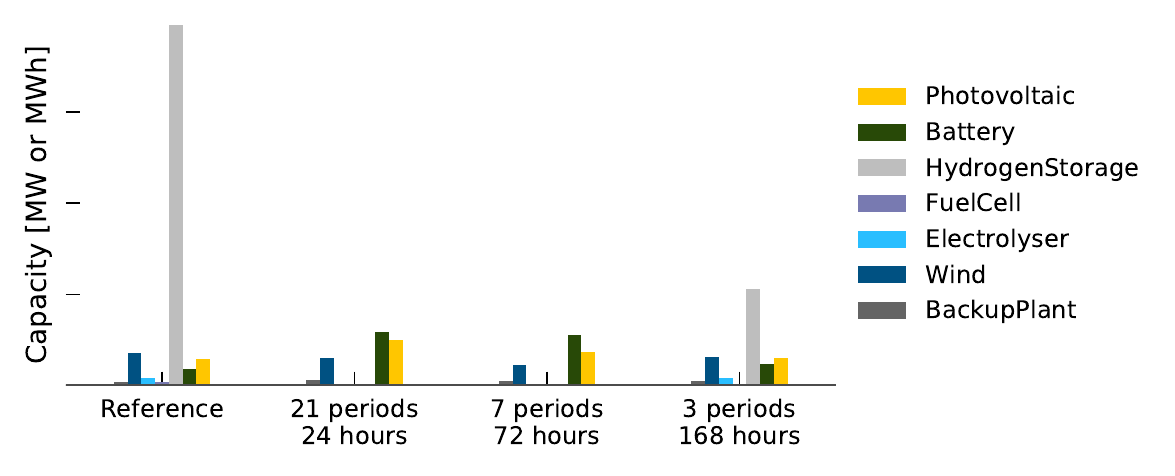}
	\caption{Scaling of the system components for an optimal island system design for the full time series and the optimal system design based on 21 typical days, seven groups of three days, and three typical weeks}
	\label{fig:Error_IntegrationMethod_capacity}
\end{figure}

The design based on the aggregation of typical days highly overestimates the  capacities of renewable energies and battery storage required. In contrast to the original system design, the hydrogen storage is not integrated at all, as is seen in Figure \ref{fig:Error_IntegrationMethod_capacity}. The same applies for the aggregation to typical groups of three days, which underestimates the overall system cost by 7\%. Only in the case of the aggregation of typical weeks does a small hydrogen storage infrastructure gets built. All in all, this system design bears the highest degree of similarity to the original system, but also underestimates the costs by 9 \%. 

In summary, a poor performance of the typical period aggregation for the highly renewable based energy system design, including storages, can be identified. The cause of this is that the sequence of typical periods is not considered in the system model, with the consequence that long term storage cannot be properly taken into account.

\section{Summary and discussion}
\label{sec:Summary}
Comparing the aggregation methods, the results show that \textit{averaging} time series based on their appearance in the year, for example the average daily profile of a month, leads to an inaccurate system design.

As alternative aggregation methods, \textit{k-mean} clustering, \textit{k-medoid} clustering and \textit{hierarchical} clustering have been compared for the creation of typical periods. In general, no aggregation method outperforms all the others for every case, but a trend is recognizable: The medoid, the most representative candidate period,  performs better than the centroid, respectively an averaged time-series. The choice of the aggregation method itself, in this case \textit{k-medoid} or \textit{hierarchical}, has a fairly small influence.

The comparison of the aggregations of different time series indicates that solar irradiation, temperature and the electricity load of a whole region can be easily represented with a small error by means of the introduced aggregation methods. An electricity load of a single building is difficult to aggregate due to high fluctuation and the small daily pattern. Furthermore, wind data is challenging: The representative periods smooth out much of the original intra-period fluctuations. All aggregations were evaluated by so-called performance indicators, which have limited analytical value.

The methods have been further applied to three different system design optimizations: A CHP-based system, a residential system based on a heat pump and photovoltaics, and an island system with a high share of renewable energy. The outcome shows that the impact of the time series aggregation is significantly reliant on the structure of the energy system itself.

The best results could be achieved for the CHP and the residential system, which are highly dependent on central supply infrastructures – either gas imports, which would also hold for other fuel imports, or a connection to the electricity grid. Therefore, a total of eight to twelve typical days was sufficient to reach cost prediction errors less than 2\% while reducing the solving duration of the optimization problem by a factor of 50 in comparison to using the time series of a full year. 

For the case of the island system, the current modeling approach leads to high estimation errors. The optimization based on the aggregated periods caused a modeling error because the typical periods are independent and cannot exchange energy. Therefore, storage technologies were not properly taken into account. An extension from typical days to typical weeks brought about a small improvement but was, as expected, still unable to account for seasonal storage.

\section{Conclusion}
\label{sec:Conclusion}
This paper investigated the effect of time series aggregation on optimal energy system design models. We showed that a time series aggregation based on clustering algorithms can significantly reduce model complexity and the required computational time. 

Comparing different aggregation methods, the results indicate that the choice of the aggregation algorithm itself has only a minor impact on the optimal system design. Nevertheless, the authors would recommend a hierarchical aggregation due to its small computational load and reproducibility. The choice of the representative periods from the aggregated clusters is of more relevance; cluster medoids, i.e., using real periods, produced more accurate results than k-mean centroids, which use averaged periods. 

A qualitative analysis using aggregated time series samples showed that the series representing entire regions can be more easily aggregated than those for single locations, i.e., less typical periods, can be sufficiently representative. Hence, aggregation methods are expected to be more valuable for spatially-aggregated systems than for models with finer resolution.

Furthermore, we illustrate that the trade-off between the depth of time series reduction and the resulting modeling error depends highly on the system configuration. While energy systems based on centralized supply resources can be well represented with a few typical days, energy systems heavily relying on storage technologies cannot be properly represented by independent typical days at all, nor by independent typical weeks. 

This leads to the conclusion that the impact of applying time series aggregation methods must be evaluated separately for each energy system model, or at least for a simplified model version. Merely assessing the aggregated time series by performance indicators can lead to significant deviation from an optimal design.

\section{Outlook}
\label{sec:Outlook}
Based on the conclusion, new modeling methods are required that take the sequence of the aggregated typical periods into account, with the objective of modeling the operation of long term storage contingencies.

Furthermore, future research should identify more cluster-relevant parameters like gradients and variance, etc., because the current clustering leads to smoothed typical periods that underestimate the variability of the original time series. 

All methods introduced are published in the Python package \href{https://github.com/FZJ-IEK3-VSA/tsam/}{tsam - Time Series Aggregation Module} and can be easily applied and extended.

\section*{Acknowledgments}
\label{s:Acknow}
This work was supported by the Helmholtz Association under the Joint Initiative "EnergySystem~2050 – A Contribution of the Research Field Energy".

\appendix
\label{sec:Appendix}
\section{System modeling}
\label{sec:systemModel}
For validation purposes, a comprehensible simple Mixed-Integer Linear Program has been chosen as the system model. 
This is defined by a network of specific technologies that are connected by energy flow variables $\dot{E}_{i,j,t}$\nomenclature[A]{$E$}{Energy flow between two components}{}{} at time step $t$. Each connection is therefore defined by an output component $i$ and input component $j$ and belongs to a connection set $(i,j) \in L$. These connections are restricted by the introduced component models.\nomenclature[B]{L}{Set of component connections}{}{} 

For a typical period, the time steps $t$ are replaced by the time steps $g$ within a single period $k$.

\subsection{Objective function}
The objective function describes the annualized cost of the supply system considered. Therefore, for each device $d$\nomenclature[C]{d}{Considered device or technology}{}{}, the annualized cost are calculated with a capital recovery factor $CRF_{d}$\nomenclature[D]{$CRF$}{Capital Recovery Factor}{}{}, which considers the Weighted Average Cost of Capital \nomenclature[D]{WACC}{Weighted Average Cost of Capital}{}{} $WACC_{d}$ and lifetime $\tau_{d}$\nomenclature[B]{$\tau$}{Lifetime}{}{} of the device in years:
\begin{align}
CRF_{d}=\frac{(1-WACC_{d})^{\tau_{d}} WACC_{d}}{ (1-WACC_{d})^{\tau_{d}}-1}
\end{align}
The capital expenditure of each component is divided into the existing related costs [eur], which only appear if the component is installed, and scale related costs [eur/kW], specific costs, which are scale dependent \cite{Lindberg2016}. For this reason, each component is modeled by a binary variable $\delta_d$ that defines whether the component exists, and a continuous variable $D_d$, which defines the installed capacity of the component. The resulting device specific annualized fixed cost can be calculated with the existing related capital expenditure ($CAPEX_{exist}$), the scaling-related capital expenditure ($CAPEX_{spec}$)\nomenclature[D]{CAPEX}{Specific capital expenditure}{}{} and fixed operational expenditure ($OPEX_{fix,d}$)\nomenclature[D]{OPEX}{Specific fix operational expenditure}{}{}  as follows 
\begin{align}
\begin{array}{c}
c_{exist,d} =  CAPEX_{exist,d} \left( CRF_{d} + OPEX_{fix,d} \right) \\
c_{spec,d} =  CAPEX_{spec,d} \left( CRF_{d} + OPEX_{fix,d} \right) \\
\end{array}
\end{align}
The costs that variate with the operation of the system  $c_{var,i,j,t}$ are related to the energy flows $\dot{E}_{i,j,t}$. Along with the scaling of the devices $D_{d}$\nomenclature[A]{$D$}{Scaling of a device}{}{}, the following objective function can be stated:
\begin{align}
\min \sum_{d} c_{exist,d} \delta_{d} + c_{spec,d} D_{d} + \sum_{(i,j) \in L} \sum_{t \in T} c_{var,i,j,t} \dot{E}_{i,j,t} \triangle t
\end{align}

\subsection{Constraints}
The device models establish the constraints of the system. They are divided into five classes: \textit{Source/Sinks, Collectors,  Transformers} and \textit{Storages}. 

The \textit{Source/Sink} class $q$\nomenclature[C]{$q$}{Index of the Source/Sink class}{}{} represents input and output flows to the system, like photovoltaic feed-in or electricity demand. It is essentially defined by a single equation:
\begin{align}
\eta_{lb,q,t} D_q \leq \sum_{(q,j) \in L } \dot{E}_{q,j,t} \leq \eta_{ub,q,t} D_q \quad \forall \quad t,q
\end{align}
where $\eta_{lb,q,t}D_q$\nomenclature[D]{$LB$}{Lower bound}{}{} \nomenclature[D]{$UB$}{Upper bound}{}{}\nomenclature[B]{$\eta$}{Efficiency}{}{}could be a certain demand that must at least be satisfied at timestep $t$, or $\eta_{ub,q,t}D_q$ could be the maximal photovoltaic feed-in per installed capacity.

The \textit{Collectors} class $n$\nomenclature[C]{$n$}{Index of the Collector class}{}{} can be seen as a hub where all input energy flows must be equivalent to all output energy flows:
\begin{align}
\sum_{(i,n) \in L } \dot{E}_{i,n,t} - \sum_{(n,j) \in L } \dot{E}_{n,j,t} = 0 \quad \forall \quad t,n
\end{align}

The \textit{Transformer} class $f$\nomenclature[C]{$f$}{Index of the Transformer class}{}{} represents devices that transform the energy from one form to another. Examples are a fuel cells or heat pumps. For their definition, the energy type (electricity, gas, etc.) $\epsilon$\nomenclature[C]{$\epsilon$}{Energy type}{}{} must be introduced. Each energy flow $\dot{E}_{i,j,t}$ has a certain energy type $\epsilon$. With the energy type specific transformation efficiency $\eta_{f,\epsilon_{in},\epsilon_{out}}$, the following equation can be stated for each energy transformation in the device:
\begin{align}
\eta_{f,\epsilon_{in},\epsilon_{out}} \sum_{(i,f) \in L,\epsilon_{in} } \dot{E}_{i,f,t} - \sum_{(f,j) \in L,\epsilon_{out} } \dot{E}_{f,j,t} = 0 \quad \forall \quad t,f
\end{align}

The \textit{Storage} class $s$\nomenclature[C]{$s$}{Index of the Storage class}{}{} is defined by an additional variable
the State of Charge $SOC_{s,t}$\nomenclature[A]{$SOC$}{State of charge}{}{} at time step $t$. We can utilize the Euler method to state for the state of charge in the next time step $SOC_{s,t+1}$:
\begin{equation}
\begin{array}{rl}
SOC_{s,t+1} = & SOC_{s,t} (1-\eta^{self}_{s} \Delta t ) \\
& + \eta_{s}^{char} \sum\limits_{ (i,s) \in L} \dot{E}_{i,s,t} \Delta t \\
& -  \frac{1}{\eta^{dis}_{s}} \sum\limits_{ (s,j) \in L } \dot{E}_{s,j,t} \Delta t   \\ 
\end{array}  \quad \forall \quad t,s 
\end{equation}
where $\dot{E}^{char}_{s,t}$ describes the charging flow with an efficiency of $\eta^{char}_{s}$ and $\dot{E}^{dis}_{s,t}$ the discharging flow with related efficiency $\eta^{dis}_{s}$. $\eta^{self}_{s}$ defines the self-discharge of the storage and $\Delta t$\nomenclature[B]{$\Delta t$}{Duration of a single time step}{}{} the step length of a single time step. The state of charge at the beginning of the considered time frame $SOC_{s,1}$ is related to that in the end of the time frame $SOC_{s,N_t+1}$.

The design variable of the storage $s$ which is described by its capacity $D_s$ limits the state of charge to the following:
\begin{equation}
SOC_{s,t} \leq D_s \quad \forall \quad t,s 
\end{equation}

The existing related variable $\delta_d$\nomenclature[A]{$\delta$}{Binary variable determining the existance of a technology}{}{} restricts the scaling-dependent device variable $D_d$ by the so called BigM-Method \cite{Bemporad1999} as follows:
\begin{equation}
\mathbf{M}  \delta_d \geq D_d
\end{equation}
The method is inspired by \citet{Stadler2014} and \citet{Lindberg2016}.

\section*{References}

\bibliography{mybibfile}

\end{document}